\documentclass{article}
\usepackage[T1]{fontenc}
\usepackage[latin9]{inputenc}
\usepackage{mathtools}
\usepackage{enumitem}
\usepackage{amsmath}
\usepackage{amsthm}
\usepackage{amssymb}
\usepackage{graphicx}
\usepackage[authoryear]{natbib}

\makeatletter
\theoremstyle{plain}
\newtheorem{thm}{\protect\theoremname}
\theoremstyle{plain}
\newtheorem{prop}[thm]{\protect\propositionname}
\theoremstyle{remark}
\newtheorem{rem}[thm]{\protect\remarkname}
\theoremstyle{plain}
\newtheorem{lem}[thm]{\protect\lemmaname}

\@ifundefined{date}{}{\date{}}
\usepackage[a4paper, includeheadfoot, margin = 2.54cm]{geometry}

\usepackage{bbm}

\usepackage{mathtools}

\usepackage{placeins}

\allowdisplaybreaks

\usepackage{hyperref}
\usepackage[capitalise,nameinlink]{cleveref}
\usepackage{autonum}

\makeatletter
\newcommand{\restore@Environment}[1]{%
  \AtBeginDocument{%
    \csletcs{#1*}{#1}%
    \csletcs{end#1*}{end#1}%
  }%
}
\forcsvlist\restore@Environment{alignat,equation,gather,multline,flalign,align}
\makeatother

\crefname{equation}{}{}

\crefalias{thm}{theorem}

\usepackage{etoolbox}

\ifcsmacro{defn}{}{\newtheorem{defn}{test}}
\ifcsmacro{thm}{}{\newtheorem{thm}{test}}
\ifcsmacro{lem}{}{\newtheorem{lem}{test}}
\ifcsmacro{prop}{}{\newtheorem{prop}{test}}
\ifcsmacro{cor}{}{\newtheorem{cor}{test}}
\ifcsmacro{assumption}{}{\newtheorem{assumption}{test}}
\ifcsmacro{example}{}{\newtheorem{example}{test}}

\let\oldthm\thm
\renewcommand{\thm}{%
\crefalias{thm}{theorem}
\oldthm
}

\let\oldlem\lem
\renewcommand{\lem}{%
\crefalias{thm}{lemma}
\oldlem
}

\let\oldprop\prop
\renewcommand{\prop}{%
\crefalias{thm}{proposition}
\oldprop
}

\let\oldcor\cor
\renewcommand{\cor}{%
\crefalias{thm}{corollary}
\oldcor
}

\let\oldexample\example
\renewcommand{\example}{%
\crefalias{thm}{example}
\oldexample
}

\let\oldassumption\assumption
\renewcommand{\assumption}{%
\crefalias{thm}{assumption}
\oldassumption
}

\usepackage{crossreftools}
\let\ORGhypersetup\hypersetup
\protected\def\hypersetup{\ORGhypersetup}
\AtBeginDocument{
\let\ref\cref
}

\setenumerate[0]{label=(\roman*)}

\usepackage{totcount}

\newtotcounter{citnum} 
\def\oldbibitem{} \let\oldbibitem=\bibitem
\def\bibitem{\stepcounter{citnum}\oldbibitem}

\makeatother

\usepackage[english]{babel}
\providecommand{\lemmaname}{Lemma}
\providecommand{\propositionname}{Proposition}
\providecommand{\remarkname}{Remark}
\providecommand{\theoremname}{Theorem}

\begin{document}
\global\long\def\N{\mathbb{N}}%
\global\long\def\Z{\mathbb{Z}}%
\global\long\def\Q{\mathbb{Q}}%
\global\long\def\R{\mathbb{R}}%
\global\long\def\C{\mathbb{C}}%
\global\long\def\P{\mathbb{P}}%
\global\long\def\E{\mathbb{E}}%
\global\long\def\Var{\operatorname*{Var}}%
\global\long\def\le{\leqslant}%
\global\long\def\ge{\geqslant}%
\global\long\def\leq{\leqslant}%
\global\long\def\geq{\geqslant}%
\global\long\def\d{\mathrm{d}}%
\global\long\def\gesssim{\gtrsim}%
\global\long\def\e{\mathrm{e}}%
\global\long\def\1{\mathbbm1}%
 
\global\long\def\eps{\varepsilon}%
\global\long\def\theta{\vartheta}%
\global\long\def\phi{\varphi}%
\global\long\def\tilde#1{\widetilde{#1}}%
\global\long\def\hat#1{\widehat{#1}}%
\global\long\def\argmin#1{\operatorname*{\mathrm{arg\,min}}_{#1}}%
\global\long\def\argmax#1{\operatorname*{\mathrm{arg\,max}}_{#1}}%
\global\long\def\arginf#1{\operatorname*{\mathrm{arg\,inf}}_{#1}}%
\global\long\def\argsup#1{\operatorname*{\mathrm{arg\,sup}}_{#1}}%
\global\long\def\diam{\operatorname*{diam}}%
\global\long\def\spann{\operatorname*{span}}%
\global\long\def\tr{\operatorname*{tr}}%
\global\long\def\supp{\operatorname*{supp}}%
\global\long\def\Bigtimes{\bigtimes}%
\global\long\def\comp{\mathsf{c}}%
\global\long\def\lebesgue{\lambda}%
\global\long\def\Re{\operatorname*{Re}}%
\global\long\def\Im{\operatorname*{Im}}%
\global\long\def\sgn{\operatorname*{sgn}}%
\global\long\def\del{\tilde{\delta}}%
\global\long\def\tg{\delta}%

\global\long\def\D{\mathcal{D}}%
\global\long\def\F{\mathcal{F}}%
\global\long\def\i{\mathrm{i}}%
\global\long\def\MS{\blacklozenge}%
\global\long\def\lap{\Delta}%

\title{Estimating a multivariate Lévy density \\
based on discrete observations }
\author{Maximilian F. Steffen\thanks{The author thanks Mathias Trabs for helpful comments. Financial support
of the DFG through project TR 1349/3-1 is gratefully acknowledged.
The author acknowledges support by the state of Baden-Württemberg
through bwHPC.}}
\date{Karlsruhe Institute of Technology}
\maketitle
\begin{abstract}
\noindent Existing results for the estimation of the Lévy measure
are mostly limited to the onedimensional setting. We apply the spectral
method to multidimensional Lévy processes in order to construct a
nonparametric estimator for the multivariate jump distribution. We
prove convergence rates for the uniform estimation error under both
a low- and a high-frequency observation regime. The method is robust
to various dependence structures. Along the way, we present a uniform
risk bound for the multivariate empirical characteristic function
and its partial derivatives. The method is illustrated with simulation
examples.
\end{abstract}
\textbf{Keywords:} Lévy processes, jump processes, multivariate density
estimation, spectral methods, 

low-frequency, high-frequency

\smallskip{}

\noindent\textbf{MSC 2020:} 62G05, 62G07, 62M15, 60G51

\section{Introduction}

Lévy processes are a staple to model continuous-time phenomena involving
jumps, for instance in physics and finance, see \citet{woyczynski2001levy}
and \citet{Cont2004}, respectively. Naturally, manifold such applications
call for multivariate processes which significantly complicates the
theoretical analysis compared to the onedimensional case. Matters
worsen as practitioners often only have time-discrete data at their
disposal which obstructs the identification of the jumps and hence
the calibration of such models. Statistical results in this setting
are typically limited to the onedimensional case or omit the estimation
of the jump distribution, despite the practical relevance. In the
present work, we study exactly this problem: estimating the jump distribution
of a multivariate Lévy process based on discrete observations.

On a theoretical level, the distribution of the Lévy process is uniquely
determined by its characteristic triplet, that is, the volatility-matrix,
the drift and the Lévy measure. The latter characterizes the jump
distribution we are interested in. From a statistical point of view,
the estimation of the Lévy measure is most challenging as we are faced
with a nonparametric problem. 

The literature commonly distinguishes the following two observation
regimes for a Lévy process observed at equidistant time points $0<\tg,2\tg,\dots,n\tg\eqqcolon T$:
Under  the low-frequency regime, $\tg$ is fixed as $n\to\infty$,
whereas $\tg\searrow0$ under the high-frequency regime. Our estimation
method is robust across sampling frequencies.

Motivated by the clearer separation between the jumps themselves and
the Gaussian component as $\tg\searrow0$ under the high-frequency
regime, threshold-based estimators have been applied extensively.
Beyond the overview given by \citet{AitSahalia2012}, \citet{Duval2021}
apply such an approach to the estimation of the Lévy measure, \citet{Gegler2010}
study the estimation of the entire Lévy triplet and \citet{Mies2020}
estimates the Blumenthal-Getoor index. However, these references are
restricted to the onedimensional case and multidimensional extensions
seem tedious due to the multitude of directions in which the process
can jump. A notable exception is the work by \citet{Buecher2013},
who estimate the tail-integrals of a multivariate Lévy process.

Under the low-freqency regime, we cannot identify the intermittent
jumps even in the absence of a Gaussian component resulting in an
ill-posed inverse problem, see \citet{Neumann2009}. A popular way
out is the spectral method, see \citet{belomestny2015}, which leverages
the relationship of the Lévy triplet with the characteristic function
of the process at any time point. Turning the observations of the
process into increments, this characteristic function is estimated
and then used to draw inference on parts of the Lévy triplet. The
method was first considered by \citet{Belomestny2006} in the context
of exponential Lévy models and has since been studied extensively,
see \citet{Belomestny2010}, \citet{gugushvili2012}, \citet{Nickl2012},
\citet{Reiss2013}, \citet{trabs2015}.

We adapt this approach to the multivariate setting by constructing
a nonparametric estimator for the Lévy density $\nu$, assuming that
it exists. Our estimator requires no knowledge of the volatility and
drift parameters and works uniformly over fully nonparametric classes
of Lévy processes with mild assumptions. In particular, Lévy processes
with infinite jump activity are allowed. The uniform rates we achieve
naturally extend those from the onedimsional case and optimality in
our setting is discussed. 

When estimating the Lévy density close to the origin, we enhance our
method with an estimator for the volatility. The estimation of the
volatility matrix itself has previously been studied, see \citet[high-frequency]{papagiannouli2020}
and \citet[low-frequency]{belomestny2018}. A related issue is the
estimation of the covariance matrix in deconvolution problems, see
\citet{belomestny2019}. However, even the proven minimax-optimal
rates of convergence are too slow as to not affect our overall rates
under the low-frequency regime. It is sufficient to estimate the trace
of the volatility matrix and we show that this can be done with a
much faster rate. With this enhancement, there is no additional loss
in the rate for the estimation of the Lévy density.

An effect emerging for multivariate processes is the possibility of
different dependence structures between the components which can be
in disagreement with the existence of a Lévy density in the form of
a Lebesgue density on the whole state space. Statistical results in
such settings are even rarer. \citet{Belomestny2011} estimates the
Lévy density for a time changed Lévy process with independent components.
We propose a quantification of the estimation error when integrating
against regular test functions under various forms of dependence structures
without modifications to our method.

The paper is organized as follows. In \ref{sec:results}, we introduce
the estimation method and state our main results along with a short
outline of the proof and the key tools used. The empirical performance
of our estimator is illustrated in a simulation examples in \ref{sec:simulation}.
The full proofs are postponed to \ref{sec:proofs}.

\section{Estimation method and main results\label{sec:results}}

We begin by introducing some notation: Throughout, an $\R^{d}$-valued
Lévy process $(L_{t})_{t\ge0}$ with Lévy measure $\nu$ is observed
in the form of of $n\in\N$ increments at equidistant time points
with time difference $\tg>0$ and overall time horizon $T\coloneqq n\tg$:
\[
Y_{k}\coloneqq L_{\tg k}-L_{\tg(k-1)},\qquad k=1,\dots,n.
\]
For $x,y\in\C^{d}$ and $p\ge1$, set $\vert x\vert_{p}\coloneqq\big(\sum_{k=1}^{d}\vert x_{k}\vert^{p}\big)^{1/p}$,
$\vert x\vert\coloneqq\vert x\vert_{2}$, $\vert x\vert_{\infty}\coloneqq\max_{k=1,\dots,d}\vert x_{k}\vert$,
$x\cdot y\coloneqq\langle x,y\rangle\coloneqq\sum_{k=1}^{d}x_{k}y_{k}$
and $x^{2}\coloneqq x\cdot x$. For a multi-index $\beta\in\N_{0}^{d}$,
we set $x^{\beta}\coloneqq\prod_{k=1}^{d}x_{k}^{\beta_{k}}$, $\vert x\vert^{\beta}\coloneqq\prod_{k=1}^{d}\vert x_{k}\vert^{\beta_{k}}$. 

If $\int\vert x\vert^{2}\,\nu(\d x)<\infty$, then the characteristic
function of $L_{t}$ is given by
\[
\phi_{t}(u)\coloneqq\E[\e^{\i\langle u,L_{t}\rangle}]=\e^{t\psi(u)}\qquad\text{with}\qquad\psi(u)\coloneqq\i\langle\gamma,u\rangle-\frac{1}{2}\langle u,\Sigma u\rangle+\int\big(\e^{\i\langle u,x\rangle}-1-\i\langle u,x\rangle\big)\,\nu(\d x)
\]
for some drift parameter $\gamma\in\R^{d}$ and some positive semidefinite
volatility matrix $\Sigma\in\R^{d\times d}$, see \citet{sato1999}.
Denoting by $\lap$ the Laplace operator, i.e.
\[
\lap g\coloneqq\sum_{k=1}^{d}\frac{\partial^{2}g}{\partial u_{k}^{2}}
\]
for a function $g\colon\R^{d}\to\C$ which is twice differentiable
in every direction, we have 
\begin{align}
\nabla\psi(u) & =\i\gamma-\Sigma u+\i\int x\big(\e^{i\langle u,x\rangle}-1\big)\,\nu(\d x),\label{eq:psigrad}\\
\lap\psi(u) & =-\tr(\Sigma)-\int\vert x\vert^{2}\e^{\i\langle u,x\rangle}\,\nu(\d x)=-\tr(\Sigma)-\F[\vert x\vert^{2}\nu](u)=\frac{\phi_{t}(u)\lap\phi_{t}(u)-(\nabla\phi_{t}(u))^{2}}{t\phi_{t}^{2}(u)},\label{eq:psilap}
\end{align}
where the integral in the first line is component-wise and $\F[\vert x\vert^{2}\nu]\coloneqq\int\e^{\i\langle\cdot,x\rangle}\vert x\vert^{2}\,\nu(\d x)$.

To motivate our estimator for $\nu$, suppose $\Sigma=0$. In view
of \ref{eq:psilap}, we then have $\nu=-\vert\cdot\vert^{-2}\F^{-1}[\lap\psi]$
and $\lap\psi$ can naturally be estimated using the empirical characteristic
function $\hat{\phi}_{\tg,n}(u)\coloneqq\frac{1}{n}\sum_{k=1}^{n}\e^{\i\langle u,Y_{k}\rangle}$
leading to
\begin{equation}
\hat{\lap\psi_{n}}(u)\coloneqq\frac{\hat{\phi}_{\tg,n}(u)\lap\hat{\phi}_{\tg,n}(u)-(\nabla\hat{\phi}_{\tg,n}(u))^{2}}{\tg\hat{\phi}_{\tg,n}(u)}\1_{\{\vert\hat{\phi}_{\tg,n}(u)\vert\ge T^{-1/2}\}}\label{eq:psilaphat}
\end{equation}
with the indicator ensuring a well-defined expression. Therefore,
granted $\nu$ has a Lévy density also denoted by $\nu$, it is reasonable
to propose the estimator 
\begin{equation}
\hat{\nu}_{h}(x)\coloneqq-\vert x\vert^{-2}\F^{-1}\big[\F K_{h}\hat{\lap\psi_{n}}\big](x),\qquad x\in\R^{d}\setminus\{0\},\label{eq:estimator}
\end{equation}
where $K$ is kernel a limited by bandwidth $h>0$ ($K_{h}\coloneqq h^{-d}K(\cdot/h)$)
satisfying for some order $p\in\N$ that for any multi-index $0\ne\beta\in\N_{0}^{d}$
with $\vert\beta\vert_{1}\le p$ we have

\begin{equation}
\begin{gathered}\int_{\R^{d}}K(x)\,\d x=1,\qquad\int x^{\beta}K(x)\,\d x=0\qquad\text{and}\qquad\supp\F K\subseteq[-1,1]^{d}\end{gathered}
.\label{eq:kernel}
\end{equation}
For $d=1$, we recover the jump density estimator by \citet{trabs2015}
to estimate quantiles of Lévy measures.

A suitable kernel can be constructed as $K\coloneqq(\F^{-1}g)/g(0)$
from an integrable even function $g\colon C^{\infty}(\R^{d})\to\R$
with support contained in $[-1,1]^{d}$, $g(0)\ne0$ and vanishing
mixed partial derivatives of order up to $p$ at $0$. For the theoretical
analysis, it will be useful to consider a kernel with product structure
$K(x)=\prod_{j=1}^{d}K^{j}(x_{j})$ for kernels $K^{j}$ on $\R$,
each with order $p$, i.e. for all $q\in\N$, $q\le p$
\[
\int_{\R}K^{j}(x_{j})\,\d x_{j},\qquad\int x_{j}^{q}K^{j}(x_{j})\,\d x_{j}=0\qquad\text{and}\qquad\supp\F K^{j}\subseteq[-1,1].
\]
Obviously, such a product kernel also fulfills \ref{eq:kernel}.

\subsection{Convergence rates}

To control the estimation error, we need to impose smoothness and
moment conditions on the Lévy density. To this end, we introduce for
a number of moments $m>0$, a regularity index $s>0$, an open subset
$U\subseteq\R^{d}$ and a universal constant $R>0$
\begin{align*}
\mathcal{C}^{s}(m,U,R) & \coloneqq\Big\{(\Sigma,\gamma,\nu)\,\Big|\,\Sigma\in\R^{d\times d}\text{ positive-semidefinite},\tr(\Sigma)\le R,\,\gamma\in\R^{d},\,\int\vert x\vert^{m}\nu(\d x)\le R\\
 & \qquad\qquad\qquad\qquad\qquad\qquad\qquad\qquad\qquad\nu\text{ has a Lebesgue density with }\Vert\vert x\vert^{2}\nu\Vert_{C^{s}(U)}\le R\Big\},
\end{align*}
where
\[
\Vert f\Vert_{C^{s}(U)}\coloneqq\sum_{\vert\beta\vert_{1}\le\lfloor s\rfloor}\sup_{x\in U}\vert f^{(\beta)}(x)\vert+\max_{\vert\beta\vert_{1}=\lfloor s\rfloor}\sup_{x,y\in U,x\ne y}\frac{\vert f^{(\beta)}(x)-f^{(\beta)}(y)\vert}{\vert x-y\vert^{s-\lfloor s\rfloor}}
\]
denotes the Hölder norm with regularity index $s$ for any function
$f$ which has derivatives $f^{(\beta)}$ of order up to $\lfloor s\rfloor\coloneqq\max\{k\in\N_{0}:k<s\}$
on $U$. $C^{s}(U)$ denotes the set of all Hölder-regular functions
on $U$ with regularity index $s>0$. Since we will require regularity
of the Lévy density in a small $\zeta$-neighborhood beyond $U$ for
a uniform rate, we set $U_{\zeta}\coloneqq\{x\in\R^{d}\mid\exists u\in U:\vert x-u\vert<\zeta\}$
for some radius $\zeta>0$. 

In view of \ref{eq:psilaphat} it is natural, that the estimation
error will also depend on the decay behavior of the characteristic
function, which in turn, is affected by the presence of a Gaussian
component. Therefore, we distinguish the following two classes of
Lévy processes. First, is the so-called mildly ill-posed case for
a decay exponent $\alpha>0$ 

\[
\mathcal{D}^{s}(\alpha,m,U,R,\zeta)\coloneqq\big\{(0,\gamma,\nu)\in\mathcal{C}^{s}(m,U_{\zeta},R)\,\big\vert\,\Vert(1+\vert\cdot\vert_{\infty})^{-\alpha}/\phi_{1}\Vert_{\infty}\le R,\,\Vert\vert x\vert\nu\Vert_{\infty}\le R\big\}.
\]
As alluded to in the introduction, a Gaussian component overclouds
the jumps in addition to the discrete observations and is therefore
treated as the severely ill-posed case for $\alpha,r,\eta>0$ 
\[
\mathcal{E}^{s}(\alpha,m,U,r,R,\zeta,\eta)\coloneqq\big\{(\Sigma,\gamma,\nu)\in\mathcal{C}^{s}(m,U_{\zeta},R)\,\big\vert\,\Vert\exp(-r\vert\cdot\vert_{\infty}^{\alpha})/\phi_{1}\Vert_{\infty}\le R,\vert x\vert^{3-\eta}\nu(x)\le R\;\forall\vert x\vert\le1\big\}.
\]
The parameters $\alpha$ and $r$ control the exponential decay of
the characteristic funtion. Note that $\Sigma\ne0$ already implies
$\alpha=2$.

In the mildly ill-posed case, the Blumenthal-Getoor index of the Lévy
process is at most $1$, whereas in the severely ill-posed case it
is at most $\big((3-\eta)\land2\big)\lor0$, where we set $a\land b\coloneqq\min\{a,b\}$
and $a\lor b\coloneqq\max\{a,b\}$ for $a,b\in\R$.

For these regularity classes, we are able to quantify the estimation
error as follows.
\begin{thm}
\label{thm:main}Let $\alpha,r,R,\zeta>0,s>1,m>4$ and let the kernel
satisfy \ref{eq:kernel} with order $p\ge s$. Let $U\subseteq\R^{d}$
be an open set which is bounded away from $0$. We have for $0<\tg\le R$,
$n\to\infty$: 
\end{thm}

\begin{enumerate}
\item[(M)] If $U$ is bounded and $h=h_{\tg,n}=(\log(T)/T)^{1/(2s+2\tg\alpha+d)}$,
then uniformly in $(\Sigma,\gamma,\nu)\in\qquad$ $\mathcal{D}^{s}(\alpha,m,U,R,\zeta)$
\[
\sup_{x^{\ast}\in U}\vert\hat{\nu}_{h}(x^{\ast})-\nu(x^{\ast})\vert=\mathcal{O}_{\P}\Big(\Big(\frac{\log T}{T}\Big)^{s/(2s+2\tg\alpha+d)}\Big).
\]
If $\tg=n^{-\eps}$ with $\eps\in(0,1)$, the choice $h=(\log(T)/T)^{1/(2s+d)}$
yields the rate $(\log(T)/T)^{s/(2s+d)}.$ 
\item[(S)] If $\vert\cdot\vert^{p+d}K\in L^{1}(\R^{d})$, $\eta>0$ and $h=h_{\tg,n}=(\log(T)/(4r\tg))^{-1/\alpha}$,
then uniformly in $(\Sigma,\gamma,\nu)\in\mathcal{E}^{s}(\alpha,m,U,r,R,\zeta,\eta)$
\[
\sup_{x^{\ast}\in U}\vert\hat{\nu}_{h}(x^{\ast})-\nu(x^{\ast})\vert=\mathcal{O}_{\P}\Big(\Big(\frac{\log T}{4r\tg}\Big)^{-s/\alpha}\Big).
\]
If $\tg=n^{-\eps}$ with $\frac{3}{2(s+d)+1}\lor\frac{\alpha}{2(s+d)+\alpha}<\eps<1$,
the choice $h=T^{-1/(2(s+d))}$ yields the rate $T^{-s/(2(s+d))}$.
\end{enumerate}
This theorem generalizes \citet[Proposition 2]{trabs2015} to the
multivariate case and additionally allows for high-frequency observations.
\ref{fig:compound_poisson3d,fig:compound_poisson_top} illustrate
a simulation example of the estimation method.

In the mildly ill-posed case, one can easily attain the same rates
without the logarithm when considering the pointwise loss. 

We first discuss the low-frequency regime: For $d=1$, our rates coincide
with the proven minimax-optimal rates in the corresponding nonparametric
deconvolution problems, see \citet{Fan1991}. In the mildly ill-posed
case with $d=1$, the pointwise variant of our rate has been shown
to be minimax-optimal under the assumption that $x\nu$ is $s$-Sobolev-regular,
see \citet{kappus2012}. In the severely ill-posed case with $d=1$
and $\alpha=\{1,2\}$, our rates coincide with the minimax-optimal
rates of \citet{Neumann2009}, who consider the integrated risk in
the estimation of $\Sigma\delta_{0}(\d x)+\vert x\vert^{2}(1+\vert x\vert^{2})^{-1}\nu(\d x)$
against test functions with Sobolev regularity $s$. This measure
has an atom in $0$ and is therefore not smooth. Hence, the regularity
in the rate comes purely from the test function. By considering $U$
bounded away from $0$, we can profit from the regularity of the Lévy
density outside the origin. We do not even suffer an additional loss
for the dimension in the rate, only in the constant. Therefore, the
above suggests its optimality.

\medskip{}

One sees that the rates improve as the time grid narrows. If this
refinement happens at an appropriate order compared to the growth
of the sample, the ill-posedness vanishes completely in the mildly
ill-posed case and the rate becomes polynomial in the severely ill-posed
case. In the mildly ill-posed case with high-frequency observations,
the rate corresponds to the minimax-optimal rate in a nonparametric
regression.

\medskip{}

It is straightforward to see from our proof, that when estimating
$\vert\cdot\vert^{2}\nu$ , we can forgo the exclusion of the origin
from $U$ while achieving the same rates in the mildly-ill-posed case.
In the severely ill-posed case, the unknown volatility of the Brownian
component of the Lévy process obstructs the observation of the small
jumps. Hence, we can benefit from a pilot estimator for $\Sigma$.
As discussed earlier, even a with minimax-optimal estimator for $\Sigma$,
we would suffer a loss in the overall rate. However, in view of \ref{eq:psilap},
it suffices to estimate the onedimensional parameter $\tr(\Sigma)$
which is easier compared to the $d\times d$-matrix $\Sigma$. Following
the spectral approach again, we propose the estimator

\[
\hat{\tr(\Sigma)}\coloneqq\hat{\tr(\Sigma)}_{h}\coloneqq-\int W_{h}(u)\hat{\lap\psi_{n}}(u)\,\d u.
\]
where $W_{h}=h^{d}W(h\cdot)$ for a bandwidth $h>0$ (correspoding
to the threshold $h^{-1}$) and a weight function $W\colon\R^{d}\to\R$
with 
\[
\int W(u)\,\d u=1\qquad\text{and}\qquad\supp W\subseteq[-1,1]^{d}.
\]
This estimator achieves a rate of $(\log T)^{-(s+d)/\alpha}$ and
is incorporated into the estimator for $\vert\cdot\vert^{2}\nu$ via
\[
\hat{\vert\cdot\vert^{2}\nu}_{h}\coloneqq-\F^{-1}\big[\F K_{h}\big(\hat{\lap\psi_{n}}+\hat{\tr(\Sigma)}_{h}\big)\big]
\]
leading to the following extension of \ref{thm:main}.
\begin{prop}
\label{prop:xs_nu}Let $\alpha,r,R,\zeta,\eta>0,1<s\in\N,m>4$ and
let the kernel satisfy \ref{eq:kernel} with order $p\ge s$ as well
as $\vert\cdot\vert^{p+d}K\in L^{1}(\R^{d})$. Assume $\Vert\F^{-1}[W(x)/x_{k}^{s}]\Vert_{L^{1}}<\infty$
for some $k$. Choosing $h=(\log(T)/(4r\tg))^{-1/\alpha}$ we have
uniformly in $(\Sigma,\gamma,\nu)\in\mathcal{E}^{s}(\alpha,m,\R^{d},r,R,\zeta,\eta)$
\[
\sup_{x^{\ast}\in\R^{d}}\big\vert\big(\hat{\vert\cdot\vert^{2}\nu}_{h}\big)(x^{\ast})-\vert x^{\ast}\vert^{2}\nu(x^{\ast})\big\vert=\mathcal{O}_{\P}\Big(\Big(\frac{\log T}{4r\tg}\Big)^{-s/\alpha}\Big).
\]
\begin{figure}
\centering{}\includegraphics[width=15cm]{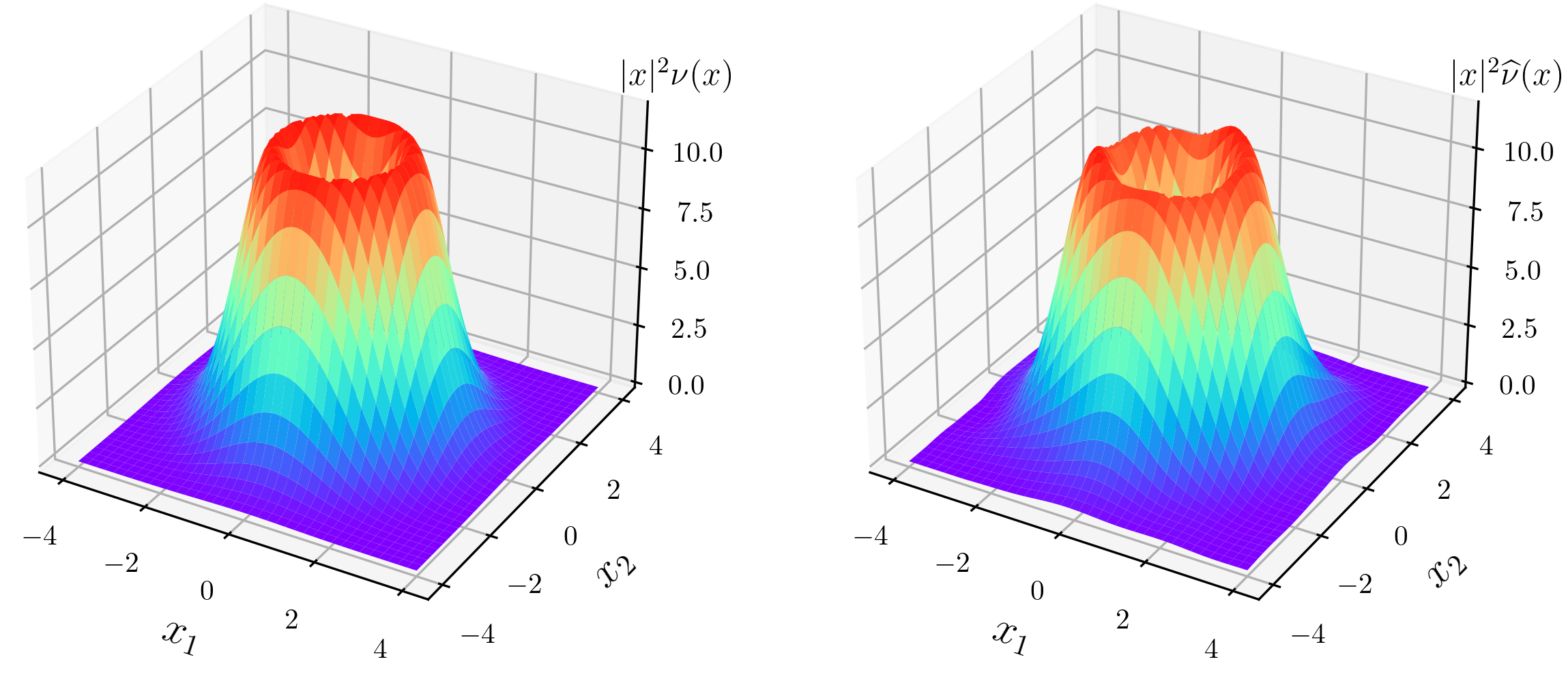}\caption{\label{fig:compound_poisson3d}3D plot of $\vert\cdot\vert^{2}\nu$
(left) and its estimate (right) for a twodimensional compound Poisson
process with Gaussian jumps.}
\end{figure}
\begin{figure}
\centering{}\includegraphics[width=15cm]{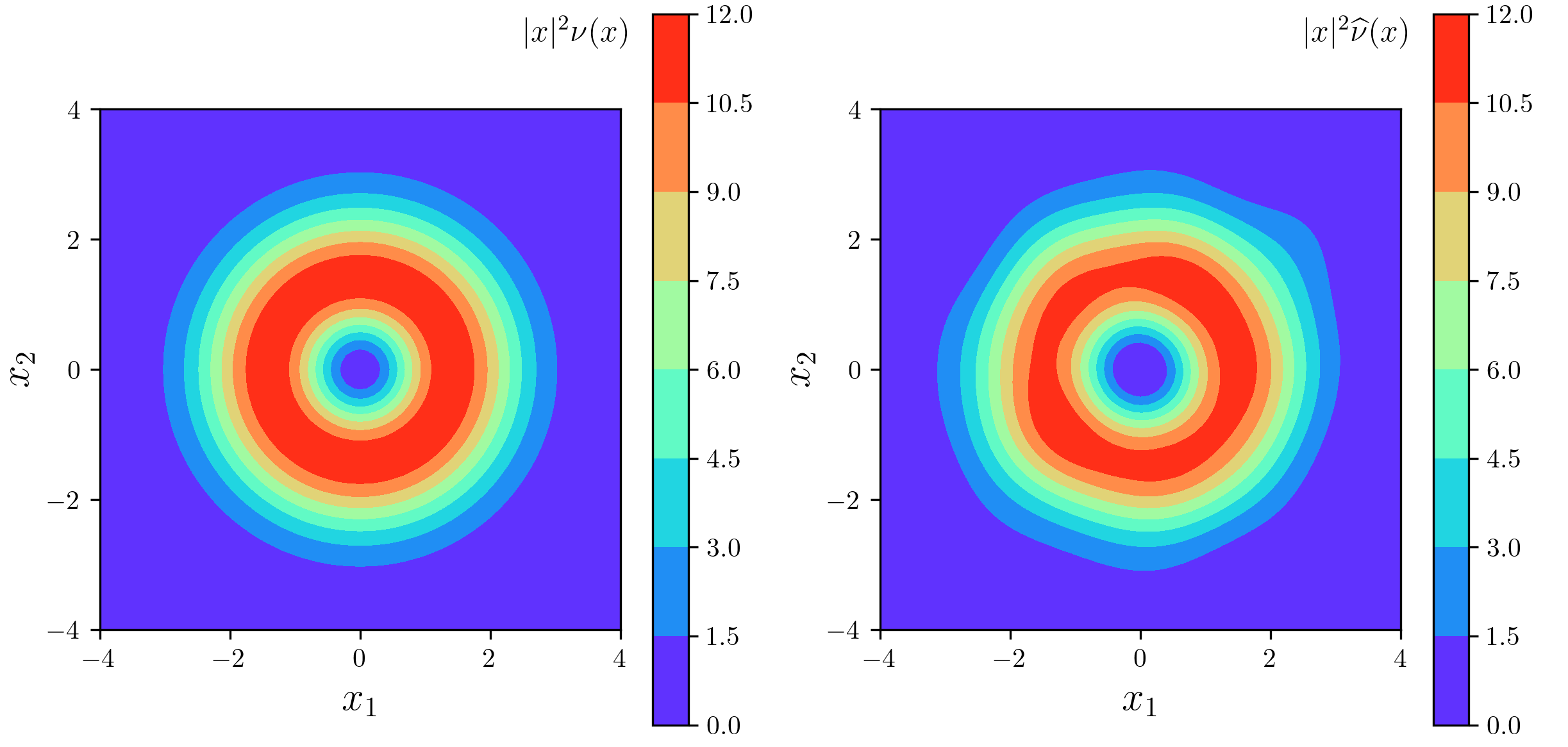}\caption{Heatmap of $\vert\cdot\vert^{2}\nu$ (left) and its estimate (right)
for a twodimensional compound Poisson process with Gaussian jumps.\label{fig:compound_poisson_top}}
\end{figure}
\end{prop}

\FloatBarrier

\subsection{Independent components\label{subsec:independent_components}}

\begin{figure}
\centering{}\includegraphics[width=15cm]{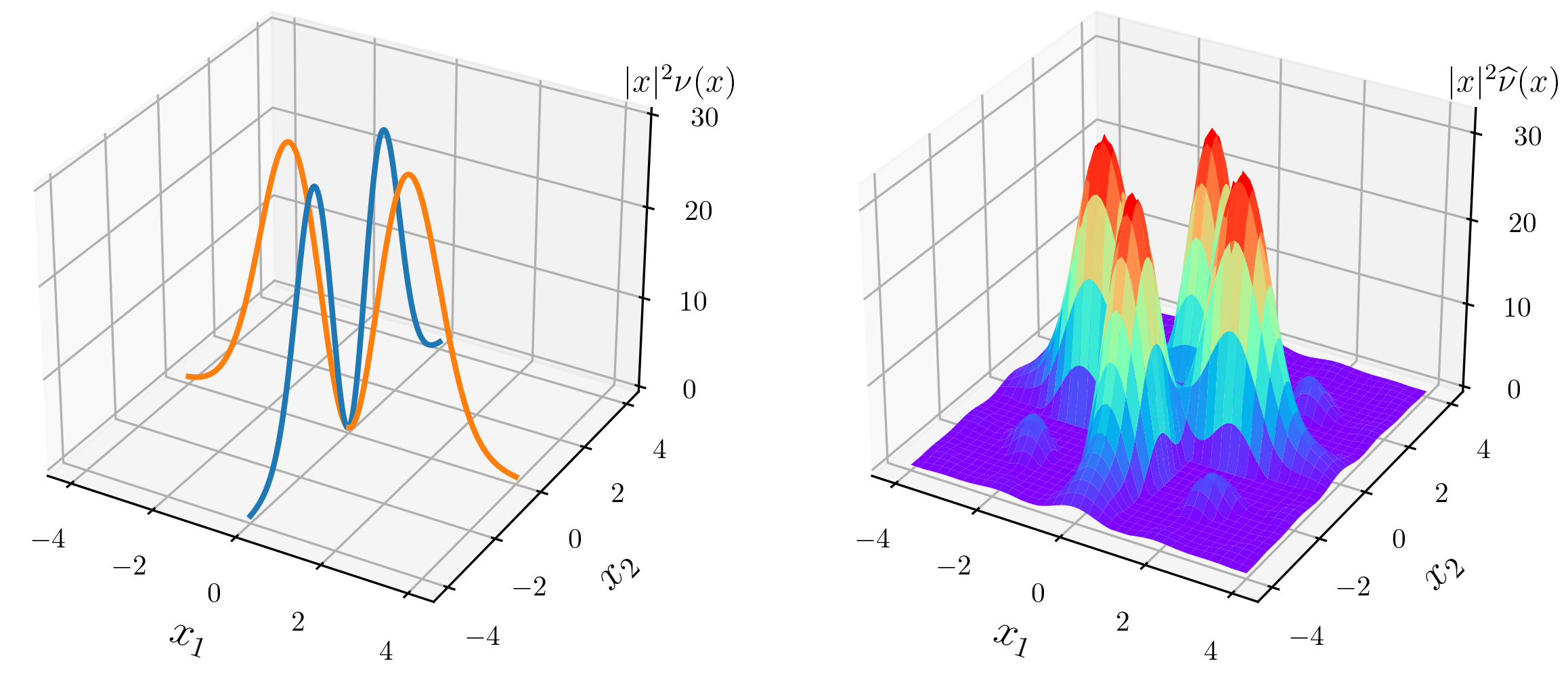}\caption{\label{fig:compound_poisson_misspec_3d}3D plot of $\vert\cdot\vert^{2}\nu$
(left) and its estimate (right) for a Lévy-process where both components
are independent compound Poisson processes with Gaussian jumps.}
\end{figure}

Compared to the onedimensional case, we need to take the dependence
structure of the components of the process into account. In particular,
our previous assumption about $\nu$ having a Lebesgue-density on
$\R^{d}$, rules out Lévy processes where all components are independent,
since the corresponding Lévy measure would only have mass on the coordinate
cross. Similarly, Lévy processes consisting of multiple mutually independent
blocks of components, where the components within the same block depend
on each other, are not covered.For the sake of notational simplicity,
we focus on the case of two equisized independent blocks: Let $d$
be even and $L=(L^{(1)},L^{(2)})$, where $L^{(1)}$ and $L^{(2)}$
are two independent Lévy processes on $\R^{d/2}$ with characteristic
triplets $(\Sigma_{1},\gamma_{1},\nu_{1})$ and $(\Sigma_{2},\gamma_{2},\nu_{2})$,
respectively. Denoting by $\delta_{0}$ the Dirac-measure in $0\in\R^{d/2}$,
it holds that
\begin{equation}
\nu(\d x)=\nu_{1}(\d x^{(1)})\otimes\delta_{0}(\d x^{(2)})+\delta_{0}(\d x^{(1)})\otimes\nu_{2}(\d x^{(2)})\qquad x=(x^{(1)},x^{(2)}),\:x^{(1)},x^{(2)}\in\R^{d/2}.\label{eq:nu_misspec}
\end{equation}
We summarize the class of such Lévy processes as 

\begin{align*}
\tilde{\mathcal{C}}(m,R) & \coloneqq\Big\{(\Sigma,\gamma,\nu)\,\Big|\,\Sigma=\Big(\!\!\!\begin{array}{cc}
\Sigma_{1} & 0\\
0 & \Sigma_{2}
\end{array}\!\!\!\Big),\tr(\Sigma)\le R,\Sigma_{1},\Sigma_{2}\in\R^{d/2\times d/2}\text{ positive-semidefinite},\gamma\in\R^{d},\\
 & \qquad\qquad\qquad\qquad\int\vert x\vert^{m}\,\nu(\d x)\le R,\nu\text{ has the form \ref{eq:nu_misspec} and \ensuremath{v_{1},\nu_{2}} have Lebesgue-densities\ensuremath{\Big\}}}
\end{align*}
for $m,R>0$. A simple example of such a Lévy measure and its estimate
are illustrated in \ref{fig:compound_poisson_misspec_3d}.

As before, we distinguish between the mildly ill-posed case with $\alpha>0$
\[
\tilde{\mathcal{D}}(\alpha,m,R)\coloneqq\big\{(0,\gamma,\nu)\in\tilde{\mathcal{C}}(m,R)\,\big\vert\,\Vert(1+\vert\cdot\vert_{\infty})^{-\alpha}/\phi_{1}\Vert_{\infty}\le R,\,\Vert\vert x_{k}\vert\nu_{k}\Vert_{\infty}\le R,k=1,2\big\}
\]
and the severely ill-posed case with $\alpha,r,\eta>0$
\[
\tilde{\mathcal{E}}(\alpha,m,r,R,\eta)\coloneqq\big\{(\Sigma,\gamma,\nu)\in\tilde{\mathcal{C}}(m,R)\,\big\vert\,\Vert\exp(-r\vert\cdot\vert_{\infty}^{\alpha})/\phi_{1}\Vert_{\infty}\le R,\vert x_{k}\vert^{3-\eta}\nu_{k}(x_{k})\le R\;\forall\vert x_{k}\vert\le1,k=1,2\big\}
\]
based on the decay behavior of the characteristic function and the
presence of a Gaussian component.

If this dependence structure were known, we could seperate the blocks
in the observations, apply our method to each block and obtain an
estimator for the overall Lévy measure. Since this is not the case,
we are left with applying our initial method. In spite of the unknown
dependence structure, we will be able to quantify the estimation error.
Due to the structure of the Lévy measure, we cannot hope for a pointwise
quantitative bound. Instead, we consider the error in a functional
sense. To this end, we introduce the following class of test functions
for $\varrho>0$ and $U\subseteq\R^{d}$ 
\[
F_{\varrho}(U,R)\coloneqq\{f\colon\R^{d}\to\R\mid f\in C^{\varrho}(\R^{d}),\Vert f\Vert_{C^{\varrho}(\R^{d})},\Vert f\Vert_{L^{1}(\R^{d})}\le R,\,\supp f\subseteq U\}.
\]
 
\begin{thm}
\label{thm:misspec} Let $\alpha,r,R>0,\varrho>1,m>4$, let the kernel
have product structure and satisfy \ref{eq:kernel} with order $p\ge\varrho$.
Then, we have for $0<\tg\le R,n\to\infty$:
\end{thm}

\begin{enumerate}
\item[(M)] If $U\subseteq\R^{d}$  is bounded and $h=(\log(T)/T)^{1/(2\varrho+2\tg\alpha+3d/2)}$,
then uniformly in $(\Sigma,\gamma,\nu)\in\tilde{\mathcal{D}}(\alpha,m,R)$
\[
\sup_{f\in F_{\varrho}(U,R)}\Big\vert\int_{U}f(x)\vert x\vert^{2}\big(\nu(\d x)-\hat{\nu}_{h}(\d x)\big)\Big\vert=\mathcal{O}_{\P}\Big(\Big(\frac{\log T}{T}\Big)^{\varrho/(2\varrho+2\tg\alpha+3d/2)}\Big).
\]
\item[(S)] If $U\subseteq\R^{d}$ is bounded away from $0$, $\vert\cdot\vert^{p+d}K\in L^{1}(\R^{d})$,
$\eta>0$ and $h=(\log(T)/(4r\tg))^{-1/\alpha}$, then uniformly in
$(\Sigma,\gamma,\nu)\in\tilde{\mathcal{E}}(\alpha,m,r,R,\eta)$ 
\[
\sup_{f\in F_{\varrho}(U,R)}\Big\vert\int_{U}f(x)\vert x\vert^{2}\big(\nu(\d x)-\hat{\nu}_{h}(\d x)\big)\Big\vert=\mathcal{O}_{\P}\Big(\Big(\frac{\log T}{4r\tg}\Big)^{-\varrho/\alpha}\Big).
\]
\end{enumerate}
Note that the regularity parameter $\varrho$ in the rates comes from
the smoothness of the test functions as compared to the smoothness
$s$ of the Lévy measure in \ref{thm:main}. In the severely ill-posed
case, the result is analogous to the well-specified. In the mildly
ill-posed case, we pay for the dependence structure with an $d/2$
in the rate. Morally, one can interpret this as the model dimension
being $3d/2$ instead of $d$. 
\begin{rem}
The product kernel is compatible with any dependence structure of
blocks, regardless of their size. For instance, if all components
of the process are independent, one gets still gets the analogous
result in severely ill-posed case. In the mildly ill-posed case, the
dimension appearing in the rate is $2d-1$ instead of $3d/2$. Comparing
the dependence structures, one finds that the two independent blocks
are an in-between case of no independent blocks and fully independent
components.
\end{rem}

\subsection{A uniform risk-bound for the characteristic function and linearization}

A key ingredient in the proofs of our preceeding results is the following
moment-uniform risk bound for the multivariate characteristic function
and its partial derivatives. It generalizes the existing results in
the univariate case (see \citealt[Theorem 1]{Kappus2010}) and the
multivariate non-uniform case (see \citealt[Proposition A.1]{belomestny2018}).
\begin{prop}
\label{lem:characteristicinequality}Let $X_{1},X_{2},\dots$ be $\R^{d}$-valued
i.i.d. random variables with characteristic function $\phi$ and empirical
characteristic function $\hat{\phi}_{n}$ such that $\E[\vert X_{1}\vert^{2\beta}\vert X_{1}\vert^{\tau}]\lesssim\rho^{\vert\beta\vert_{1}\land1}$
and $\E[\vert X_{1}\vert^{2\beta}]\lesssim\rho^{\vert\beta\vert_{1}\wedge1}$
for some multi-index $\beta\in\N_{0}^{d}$ and $\tau,\rho>0$. For
the inverse softplus-type weight function $w(u)=\log(\e+\vert u\vert)^{-(1+\chi)/2}$
with $\chi>0$, we have 
\[
\E\big[\big\Vert w(u)(\hat{\phi}_{n}-\phi)^{(\beta)}(u)\big\Vert_{\infty}\big]\lesssim\rho^{(\vert\beta\vert_{1}\land1)/2}n^{-1/2}.
\]
\end{prop}

As a direct consequence of \ref{lem:characteristicinequality}, the
indicator in the definition \ref{eq:psilaphat} equals one on the
support of $\F K_{h}$, with probability converging to one for the
bandwidths we consider.

To prove our rates for $\hat{\nu}_{h}$, \ref{eq:psilap} lets us
decompose the error 
\begin{align}
\hat{\nu}_{h} & (x^{\ast})-\nu(x^{\ast})\nonumber \\
 & \negthickspace\negthickspace\negthickspace\negthickspace=\vert x^{\ast}\vert^{-2}\big(\big(K_{h}\ast(\vert\cdot\vert^{2}\nu)-\vert\cdot\vert^{2}\nu\big)(x^{\ast})-\F^{-1}\big[\F K_{h}\big(\hat{\lap\psi_{n}}-\lap\psi\big)\big](x^{\ast})+\tr(\Sigma)K_{h}(x^{\ast})\big)\nonumber \\
 & \negthickspace\negthickspace\negthickspace\negthickspace=\vert x^{\ast}\vert^{-2}\big(\underbrace{\big(K_{h}\ast(\vert\cdot\vert^{2}\nu)-\vert\cdot\vert^{2}\nu\big)(x^{\ast})}_{\eqqcolon B^{\nu}(x^{\ast})}-\underbrace{\frac{1}{\tg}\F^{-1}\big[\F K_{h}\lap\big((\hat{\phi}_{\tg,n}-\phi_{\tg})/\phi_{\tg}\big)\big](x^{\ast})}_{\eqqcolon L_{\tg,n}^{\nu}(x^{\ast})}+R_{\tg,n}+\tr(\Sigma)K_{h}(x^{\ast})\big)\label{eq:decomposition}
\end{align}
into a bias term $B^{\nu}$, the linearized stochastic error $L_{\tg,n}^{\nu}$,
the error $\tr(\Sigma)K_{h}$ due to the volatility and a remainder
term $R_{\tg,n}$. \ref{lem:characteristicinequality} applied to
the increments of the Lévy process leads to the following linearization.
\begin{lem}
\label{lem:linearization} Let $\int\vert x\vert^{4+\tau}\,\nu(\d x)\le R$
for some $\tau>0$. If $n^{-1/2}(\log h^{-1})^{(1+\chi)/2}\Vert\phi_{\tg}^{-1}\Vert_{L^{\infty}(I_{h})}\to0$
as $n\to\infty$ for $h\in(0,1),\chi>0$, it holds
\begin{align}
\sup_{\vert u\vert_{\infty}\le h^{-1}}\big\vert\hat{\lap\psi_{n}}(u) & -\lap\psi(u)-\tg^{-1}\lap\big((\hat{\phi}_{\tg,n}-\phi_{\tg})/\phi_{\tg}\big)(u)\big\vert=\mathcal{O}_{\P}(a_{n}),\qquad\text{where}\nonumber \\
a_{n} & \coloneqq n^{-1}(\log h^{-1})^{1+\chi}\Vert\phi_{\tg}^{-1}\Vert_{L^{\infty}(I_{h})}^{2}\tg^{-1/2}\big(\tg\big\Vert\vert\nabla\psi\vert\big\Vert_{L^{\infty}(I_{h})}+\tg^{3/2}\big\Vert\vert\nabla\psi\vert\big\Vert_{L^{\infty}(I_{h})}^{2}+1\big).\label{eq:a_n}
\end{align}
\end{lem}

As a direct consequence, the remainder term is of the order 
\begin{equation}
\vert R_{\tg,n}\vert=\mathcal{O}_{\P}\big(h^{-d}a_{n}\big).\label{eq:remainderorder}
\end{equation}
After treating the four terms in \ref{eq:decomposition}, the asserted
rates follow from our bandwidth choices.

\section{Simulation examples\label{sec:simulation}}

 We demonstrate the estimation of the Lévy density for $d=2$ with
three examples: a compound Poisson process, a variance gamma process
and two independent compound Poisson processes. 

A challenge is to find examples of multivariate Lévy processes for
which paths can be simulated and the true Lévy measure is accessible
(at least numerically). To allow for plotable results, we consider
the case $d=2$ and compensate for the possible singularity of the
Lévy density at the origin, i.e. we plot $\vert\cdot\vert^{2}\nu$
and its estimate. Throughout, we use the flat-top-kernel $K$, see
\citet{McMurry2004}, as defined by its Fourier transform
\[
\F K(u)\coloneqq\begin{cases}
1, & \vert u\vert\le c,\\
\exp\Big(-\frac{b\exp(-b/(\vert u\vert-c)^{2}}{(\vert u\vert-1)^{2}}\Big) & c<\vert u\vert<1,\\
0, & \vert u\vert\ge1,
\end{cases}
\]
whose decay behaviour is controlled by $b>0$ and $0<c<1$. In our
simulations, $b=1,c=1/50$ deliver stable results. While a product
kernel is convenient for theoretical reasons in \ref{subsec:independent_components},
it did not seem necessary in practice. Throughout, we simulate increments
of the processes with a time difference of $\tg=0.001$ and fix the
bandwidth at $h=4T^{-1/2}$. To conquer this ill-posed problem, we
use large samples of $n=500,000$ increments. From the definition
\ref{eq:estimator} of the estimator, it is not guaranteed that $\hat{\nu}\ge0$
and for numerical reasons even $\hat{\nu}\in\C\setminus\R$ is possible
in practice. Therefore, we consider the estimator $\Re(\hat{\nu})\lor0$
in our simulations.

\medskip{}

The most straightforward example under consideration is the compound
Poisson process with intensity $\lambda=100$ and twodimensional standard-Gaussian
jumps. In this case, the Lévy density is just the standard normal
density, rescaled with the intensity $\lambda$. \ref{fig:compound_poisson3d}
illustrates that the method captures the overall shape of the density.
The heatmap in \ref{fig:compound_poisson_top} provides a more detailed
view especially around the origin. We observe that the decay for $\vert x\vert\to\infty$
and $\vert x\vert\searrow0$ is well-estimated, with slight problems
only arising on an annulus around the origin.\medskip{}

\begin{figure}
\centering{}\includegraphics[width=15cm]{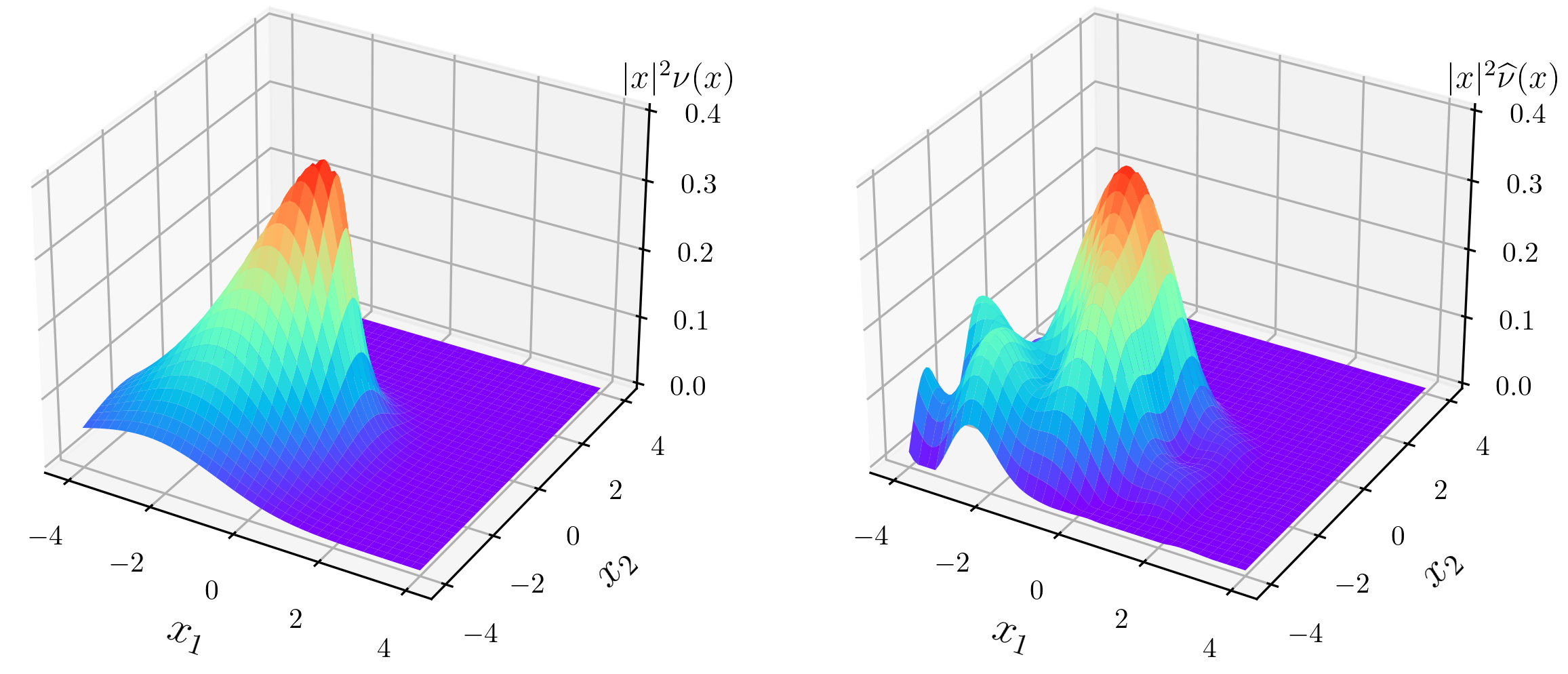}\caption{3D plot of $\vert\cdot\vert^{2}\nu$ (left) and its estimate (right)
for a twodimensional variance gamma process.\label{fig:VG_3d}}
\end{figure}

\begin{figure}
\centering{}\includegraphics[width=15cm]{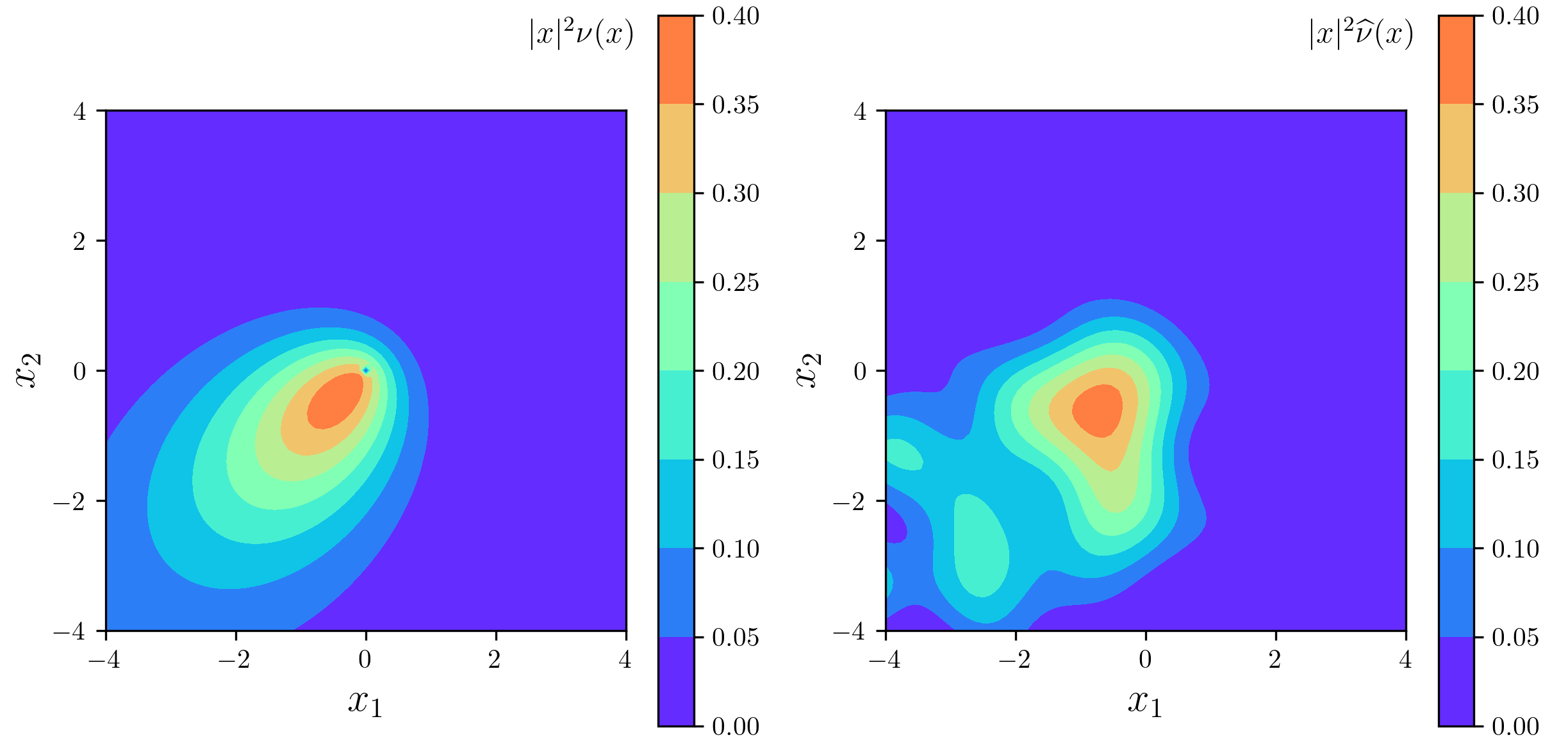}\caption{Heatmap of $\vert\cdot\vert^{2}\nu$ (left) and its estimate (right)
for a twodimensional variance gamma process.\label{fig:VG_top}}
\end{figure}

A practical way to construct easy-to-simulate multivariate Lévy processes
is to subordinate multivariate Brownian motion. In particular, we
use a gamma process with variance $\kappa=1$ to subordinate a twodimensional
standard Brownian motion. To access the Lévy measure of the resulting
variance gamma process, we approximate the theoretical expression
from \citet[Theorem 4.2]{Cont2004} numerically. The results are again
illustrated in a 3D plot (\ref{fig:VG_3d}) and as a heatmap ((\ref{fig:VG_top})).
In this example, the estimator suffers from oscillations around the
true density which are to be expected from spectral-based methods.

\medskip{}

To demonstrate the method under the depencence structure discussed
in (\ref{subsec:independent_components}), we consider a Lévy process
comprised of two independent compound Poisson processes, each with
intensity $\lambda=100$ and onedimensional standard-Gaussian jumps.
In contrast to the twodimensional compound Poisson process at the
considered at the beginnung of this section, the jumps in both components
are driven by independent Poisson processes. The corresponding Lévy
measure takes the form \ref{eq:nu_misspec}, where $\nu_{1}$ and
$\nu_{2}$ are onedimensional standard-Gaussian densities, rescaled
with $\lambda$, as illustrated on the left hand side of \ref{fig:compound_poisson_misspec_3d}.
It is important to emphasize, that the blue and the orange line represent
the Lebesgue-densities of both components on $\R$, not $\R^{2}$.
The right hand side of the aforementioned figure reveals a strong
performance of the estimator on the coordinate cross. Around the axes,
we observe a smearing effect due to the singularity of the true Lévy
measure on the coordinate cross before the estimate drops off as we
move away.

\section{Proofs\label{sec:proofs}}

Throughout, set $I_{h}\coloneqq[-h^{-1},h^{-1}]^{d}$ for $h>0$.
Note that $\Sigma,\nu,\lap\psi$ and $\hat{\lap\psi_{n}}$ do not
change if we consider increments based on the Lévy process $(L_{t}-t\gamma_{0})_{t\ge0}$
for some $\gamma_{0}\in\R^{d}$. Hence, no generality is lost if we
choose $\gamma_{0}$ such that in the mildly ill-posed case
\begin{equation}
\nabla\psi=\i\F[x\nu]\qquad\qquad\text{and}\qquad\qquad\lap\psi=-\F[\vert x\vert^{2}\nu]\label{eq:psigradwlog}
\end{equation}
and in the severely ill-posed case $\gamma=0$, see \citet[Lemma 12]{Nickl2016}
for a similar argument in the onedimensional case.

Further, due to the infinite divisibility of $\nu$, the decay behavior
of $\phi_{1}$ governs that of $\phi_{\tg}$. In particular, we have
for $0<\tg\le R$
\[
\Vert(1+\vert\cdot\vert_{\infty})^{-\tg\alpha}/\phi_{\tg}\Vert_{\infty}\le(1\vee R)^{R}\qquad\text{and}\qquad\Vert\exp(-r\tg\vert\cdot\vert_{\infty}^{\alpha})/\phi_{\tg}\Vert_{\infty}\le(1\vee R)^{R}
\]
in the mildly and the severely ill-posed case, respectively.

\subsection{Proof of \ref{thm:main}}

We extend the proof strategy by \citet{trabs2015} to accomodate for
the multivariate setting. To allow for the application to high-frequency
observations, we carefully keep track of $\tg$ throughout. Subsequently,
we will analyze the four terms in \ref{eq:decomposition}.

\subsubsection{Controlling the linearized stochastic error\label{subsec:stochasticerror}}

To control 
\[
L_{\tg,n}^{\nu}\coloneqq\F^{-1}\big[\F K_{h}\tg^{-1}\lap\big((\hat{\phi}_{\tg,n}-\phi_{\tg})/\phi_{\tg}\big)\big]
\]
we need to get a grip on the partial derivatives of $\hat{\phi}_{\tg,n}-\phi_{\tg}$
in the Laplacian of $(\hat{\phi}_{\tg,n}-\phi_{\tg})/\phi_{\tg}$.
In particular, we will show that

\begin{equation}
\sup_{u\in\R^{d}}\E\Big[\Big\vert\frac{\partial^{l}}{\partial u_{k}^{l}}(\hat{\phi}_{\tg,n}-\phi_{\tg})(u)\Big\vert\Big]\le\sup_{u\in\R^{d}}\E\Big[\Big\vert\frac{\partial^{l}}{\partial u_{k}^{l}}(\hat{\phi}_{\tg,n}-\phi_{\tg})(u)\Big\vert^{2}\Big]^{1/2}\stackrel{!}{\lesssim}n^{-1/2}\tg^{(l\land1)/2},\qquad l=0,1,2.\label{eq:moments}
\end{equation}
Since 
\[
\E\Big[\Big\vert\frac{\partial^{l}}{\partial u_{k}^{l}}(\hat{\phi}_{\tg,n}-\phi_{\tg})(u)\Big\vert\Big]\le n^{-1}\E[Y_{1,m}^{2l}]=n^{-1}\Big\vert\frac{\partial^{2l}}{\partial u_{k}^{2l}}\phi_{\tg}(0)\Big\vert\qquad\forall u\in\R^{d},
\]
where $Y_{1,k}$ denotes the $k$-th entry of $Y_{1}$, the case $l=0$
is obvious and for $l=1,2$ it remains to show
\[
\Big\vert\frac{\partial^{2l}}{\partial u_{k}^{2l}}\phi_{\tg}(0)\Big\vert\lesssim\tg.
\]
In the mildly ill-posed case 
\[
\Big\vert\frac{\partial}{\partial u_{k}}\psi(0)\Big\vert\lesssim\int\vert x_{k}\vert\,\nu(\d x)\le\int\vert x\vert\,\nu(\d x)\lesssim1
\]
and in the severely ill-posed case 
\[
\frac{\partial}{\partial u_{k}}\psi(0)=0.
\]
The product rule for higher order derivatives yields 
\begin{align*}
\Big\vert\frac{\partial^{2}\phi_{\tg}}{\partial u_{k}^{2}}(0)\Big\vert & =\Big\vert\tg\phi_{\tg}(u)\Big(\tg\Big(\frac{\partial}{\partial u_{k}}\psi(u)\Big)^{2}+\frac{\partial^{2}}{\partial u_{k}^{2}}\psi(u)\Big)\Big\vert_{u=0}\lesssim\tg\Big(1+\int\vert x\vert^{2}\,\nu(\d x)\Big)\lesssim\tg\qquad\text{and}\\
\Big\vert\frac{\partial^{4}\phi_{\tg}}{\partial u_{k}^{4}}(0)\Big\vert & =\Big\vert\frac{\partial^{2}}{\partial u_{k}^{2}}\Big(\frac{\partial^{2}\phi_{\tg}}{\partial u_{k}^{2}}(u)\Big)(0)\Big\vert\lesssim\tg\sum_{j=0}^{2}\binom{2}{j}\vert\E[Y_{1,k}^{2-j}]\vert\Big(\Big\vert\frac{\partial^{j}}{\partial u_{k}^{j}}\Big(\frac{\partial\psi}{\partial u_{k}}(u)\Big)^{2}\Big\vert_{u=0}+\Big\vert\frac{\partial^{j+2}\psi}{\partial u_{k}^{j+2}}(0)\Big\vert\Big)\lesssim\tg,
\end{align*}
where all emerging partial derivatives can again be absolutely and
uniformly bounded using our assumptions on $\nu$. 

To simplify the notation, set $m_{\tg,h}\coloneqq\F K_{h}/\phi_{\tg}$
and recall $x\cdot y\coloneqq\sum_{k=1}^{d}x_{k}y_{k}$ for $x,y\in\C^{d}$. 

For the severely ill-posed case, we have $\vert\lap\psi(u)\vert\lesssim1$
and that $\vert\e^{\i\langle u,x\rangle}-1\vert\le\vert x\vert\vert u\vert$
implies $\vert\nabla\psi(u)\vert\lesssim\vert u\vert$. Together with
\ref{eq:moments}, we obtain
\begin{align}
\E\Big[\sup_{x^{\ast}\in\R^{d}}\vert L_{\tg,n}^{\nu}(x^{\ast})\vert\Big] & \le\tg^{-1}\E\big[\big\Vert\mathcal{F}^{-1}[m_{\tg,h}\lap(\hat{\phi}_{\tg,n}-\phi_{\tg})]\big\Vert_{\infty}\big]+2\E\big[\big\Vert\F^{-1}[m_{\tg,h}\nabla(\hat{\phi}_{\tg,n}-\phi_{\tg})\cdot\nabla\psi]\big\Vert_{\infty}\big]\nonumber \\
 & \qquad\qquad\qquad+\E\big[\big\Vert\F^{-1}[m_{\tg,h}(\hat{\phi}_{\tg,n}-\phi_{\tg})(\tg(\nabla\psi)^{2}-\lap\psi)]\big\Vert_{\infty}\big]\nonumber \\
 & \lesssim(2\pi)^{-d}\int_{I_{h}}\Big(\tg^{-1}\E\big[\big\vert\lap(\hat{\phi}_{\tg,n}-\phi_{\tg})(u)\big\vert\big]+\E\big[\big\vert\nabla(\hat{\phi}_{\tg,n}-\phi_{\tg})(u)\cdot\nabla\psi(u)\big\vert\big]\nonumber \\
 & \qquad\qquad\qquad+\E\big[\big\vert(\hat{\phi}_{\tg,n}-\phi_{\tg})(u)\big\vert\big]\big(\vert\lap\psi(u)\vert+\tg\vert\nabla\psi(u)\vert^{2}\big)\Big)\exp(r\tg\vert u\vert_{\infty}^{\alpha})\,\d u\nonumber \\
 & \lesssim\pi^{-d}T^{-1/2}\int_{I_{h}}\big(1+\tg\vert u\vert+\tg^{1/2}+\tg^{3/2}\vert u\vert^{2}\big)\exp(r\tg\vert u\vert_{\infty}^{\alpha})\,\d u\nonumber \\
 & \lesssim T^{-1/2}\big(h^{-d}+\tg h^{-d-1}+\tg^{3/2}h^{-d-2}\big)\exp(r\tg h^{-\alpha}),\label{eq:severely_linearized}
\end{align}
which will be dominated by the bias.

In the mildly ill-posed case, the stochastic error needs to be decomposed
further into the main stochastic error
\begin{align}
M_{\tg,n}^{\nu} & \coloneqq-\frac{1}{T}\sum_{k=1}^{n}\F^{-1}\Big[m_{\tg,h}\big(\vert Y_{k}\vert^{2}\e^{\i\langle u,Y_{k}\rangle}-\E\big[\vert Y_{k}\vert^{2}\e^{\i\langle u,Y_{k}\rangle}\big]\big)\Big]\qquad\text{and}\nonumber \\
M_{\tg,n}^{\nu}-L_{\tg,n}^{\nu} & \hspace{0.2em}=2\F^{-1}\big[m_{\tg,h}\nabla(\hat{\phi}_{\tg,n}-\phi_{\tg})\cdot\nabla\psi\big]+\F^{-1}\big[m_{\tg,h}(\hat{\phi}_{\tg,n}-\phi_{\tg})\big(\lap\psi-\tg(\nabla\psi)^{2}\big)\big].\label{eq:diffstoch}
\end{align}
To control the difference \ref{eq:diffstoch}, note that $\Vert\vert x\vert\nu\Vert_{\infty}\le R$
and $\Vert\vert x\vert^{m}\nu\Vert_{L^{1}}\le R$ imply $\Vert\vert x\vert\nu\Vert_{L^{1}}$,
$\Vert\vert x\vert\nu\Vert_{L^{2}}$, $\Vert\vert x\vert^{2}\nu\Vert_{L^{2}}\lesssim1$.
Further, the support of $\F K$ and the decay behavior of $\phi_{\tg}$
ensure
\begin{equation}
\Vert m_{\tg,h}\Vert_{L^{2}}^{2}\lesssim\int\vert\F K(hu)\vert^{2}(1+\vert u\vert)^{2\tg\alpha}\,\d u\lesssim(1+h^{-1})^{2\tg\alpha}h^{-d}\lesssim h^{-2\tg\alpha-d}.\label{eq:m_order}
\end{equation}
Hence, \ref{eq:moments} and the Cauchy-Schwarz inequality together
with \ref{eq:psigradwlog} and the Plancherel theorem and \ref{eq:psigradwlog}
lead to 
\begin{align}
\E\Big[\sup_{x^{\ast}\in U}\big\vert M_{\tg,n}^{\nu}(x^{\ast})-L_{\tg,n}^{\nu}(x^{\ast})\big\vert\Big] & \le(2\pi)^{-d}\big(2\E\big[\big\Vert m_{\tg,h}\nabla(\hat{\phi}_{\tg,n}-\phi_{\tg})\cdot\nabla\psi\big\Vert_{L^{1}}\big]\nonumber \\
 & \qquad\qquad\qquad+\E\big[\big\Vert m_{\tg,h}(\hat{\phi}_{\tg,n}-\phi_{\tg})\big(\lap\psi-\tg(\nabla\psi)^{2}\big)\big\Vert_{L^{1}}\big]\big)\nonumber \\
 & \lesssim n^{-1/2}\Vert m_{\tg,h}\Vert_{L^{2}}\big(\tg^{1/2}\Vert\vert x\vert\nu\Vert_{L^{2}}+\Vert\vert x\vert^{2}\nu\Vert_{L^{2}}+\tg d^{2}\Vert\vert x\vert\nu\Vert_{L^{2}}\Vert\vert x\vert\nu\Vert_{L^{1}}\big)\nonumber \\
 & \lesssim n^{-1/2}h^{-\tg\alpha-d/2}.\label{eq:plancherel}
\end{align}

Being the sum of centered i.i.d. random variables, the main stochastic
error for fixed $x$ is controlled by Bernstein's inequality as summarized
in the following lemma.
\begin{lem}
\label{lem:bernsteinapplication}Let $\alpha,R,\zeta>0,m>4$ and $x\in\R^{d}$
and let the kernel satisfy \ref{eq:kernel} for $p\ge1$. If $(\Sigma,\gamma,\nu)\in\mathcal{D}^{s}(\alpha,m,U,R,\zeta)$,
then there exists some constant $c>0$ depending only on $R,\alpha$
and $d$ such that for any $\kappa_{0}>0$ and any $n\in\N,0<\tg\le R,h\in(0,1)$
\[
\P\big(\vert M_{\tg,n}^{\nu}(x)\vert\ge\kappa_{0}T^{-1/2}h^{-\tg\alpha-d/2}\big)\le2\exp\Big(-\frac{c\kappa_{0}^{2}}{(1+\vert x\vert^{3})\big(1+\kappa_{0}(h^{d}T)^{-1/2}\big)}\Big).
\]
\end{lem}

To establish a uniform bound for $x^{\ast}\in U$, a union bound extends
this lemma to a discretization of the bounded set $U$ and Lipschitz
continuity of $x\mapsto M_{\tg,n}^{\nu}(x)$ allows us to control
the discretization error. In particular, a standard covering argument
yields a discretization $x_{1},\dots,x_{N_{n}}\in\R^{d}$ of $U$
such that $\sup_{x^{\ast}\in U}\min_{l=1,\dots,N_{n}}\vert x^{\ast}-x_{l}\vert\le T^{-2}$,
$N_{n}\lesssim T^{2d}$ and $\max_{l=1,\dots,N_{n}}\vert x_{l}\vert\le C$
with some $C>0$ indepedent of $n$. Since
\begin{align*}
M_{\tg,n}^{\nu} & =K_{h}\ast g\qquad\text{with}\qquad g\coloneqq\tg^{-1}\F^{-1}\big[\1_{I_{h}}\phi_{\tg}^{-1}\lap(\hat{\phi}_{\tg,n}-\phi_{\tg})\big],
\end{align*}
the fundamental theorem of calculus together with the order of the
kernel ensures the Lipschitz continuity of $M_{\tg,n}^{\nu}$ via
\begin{align*}
\vert M_{\tg,n}^{\nu}(x)-M_{\tg,n}^{\nu}(y)\vert=\Big\vert\int_{0}^{1}(x-y)\cdot\nabla(K_{h}\ast g)(y+\tau(x-y))\,\d\tau\Big\vert & \le\vert x-y\vert\E[\Vert g\Vert_{\infty}]\big\Vert\vert\nabla K_{h}\vert_{1}\big\Vert_{L^{1}}\\
 & \lesssim\vert x-y\vert h^{-1}\E[\Vert g\Vert_{\infty}].
\end{align*}
Therefore, the discretization error is upper bounded by
\[
\E\Big[\sup_{x^{\ast}\in U}\min_{l=1,\dots,N_{n}}\vert M_{\tg,n}^{\nu}(x^{\ast})-M_{\tg,n}^{\nu}(x_{l})\vert\Big]\lesssim T^{-5/2}h^{-1}\int_{I_{h}}\vert\phi_{\tg}^{-1}(u)\vert\,\d u\lesssim T^{-5/2}h^{-\tg\alpha-d-1}.
\]
Combining the above with Markov's inequality yields for any $\kappa_{0}$
such that $2d<\frac{c}{6}\kappa_{0}^{2}/(1+C^{3})$ with $c$ from
\ref{lem:bernsteinapplication} and $T$ with $\kappa_{0}^{2}\log(T)/(Th^{d})\le1$
\begin{align}
\P\Big(\sup_{x^{\ast}\in U}\vert M_{\tg,n}^{\nu} & (x^{\ast})\vert>\kappa_{0}\Big(\frac{\log T}{T}\Big)^{1/2}h^{-\tg\alpha-d/2}\Big)\nonumber \\
 & \le\P\Big(\max_{l=1,\dots,N_{n}}\vert M_{\tg,n}^{\nu}(x_{l})\vert+\sup_{x^{\ast}\in U}\min_{l=1,\dots,N_{n}}\vert M_{\tg,n}^{\nu}(x^{\ast})-M_{\tg,n}^{\nu}(x_{l})\vert>\kappa_{0}\Big(\frac{\log T}{T}\Big)^{1/2}h^{-\tg\alpha-d/2}\Big)\nonumber \\
 & \le\frac{2}{\kappa_{0}}\Big(\frac{\log T}{T}\Big)^{-1/2}h^{\tg\alpha+d/2}\E\Big[\sup_{x^{\ast}\in U}\min_{l=1,\dots,N_{n}}\vert M_{\tg,n}^{\nu}(x^{\ast})-M_{\tg,n}^{\nu}(x_{l})\vert\Big]\nonumber \\
 & \qquad\qquad+2N_{n}\exp\Big(-\frac{c\kappa_{0}^{2}\log T}{2(1+C^{3})\big(2+\kappa_{0}(\log(T)/(Th^{d}))^{1/2}\big)}\Big)\nonumber \\
 & \lesssim\frac{2}{\kappa_{0}}\Big(\frac{\log T}{T}\Big)^{-1/2}h^{-d/2-1}T^{-5/2}+2\exp\Big(\Big(2d-\frac{c\kappa_{0}^{2}}{6(1+C^{3})}\Big)\log T\Big).\label{eq:bernsteinunion}
\end{align}
The second term obviously converges to $0$ as $T\to\infty$. For
the first term, $3d/2\ge d/2+1$ implies 
\[
\frac{2}{\kappa_{0}}\Big(\frac{\log T}{T}\Big)^{-1/2}h^{-d/2-1}T^{-5/2}\le\frac{2}{\kappa_{0}}\Big(\frac{\log T}{T}\Big)^{-1/2}h^{-3d/2}T^{-5/2}=\frac{2}{\kappa_{0}}\Big(\frac{\log T}{Th^{d}}\Big)^{3/2}T^{-1/2}(\log T)^{-2}
\]
and the right hand side converges to $0$ by our choice of bandwidth. 

Overall, \ref{eq:bernsteinunion,eq:plancherel} show
\[
\vert L_{\tg,n}^{\nu}\vert=\mathcal{O}_{\P}\Big(\Big(\frac{\log T}{T}\Big)^{1/2}h^{-\tg\alpha-d/2}\Big),
\]
which our bandwidths balance with the bias.

\subsubsection{Controlling the error owed to the volatility\label{subsec:volatility}}

We now consider the last term in \ref{eq:decomposition}. The mildly
ill-posed case is trivial since $\Sigma=0$. Turning to the severely
ill-posed case, we first aim to bound $\vert x\vert^{p+d}\vert K(x)\vert$.
To this end, consider 
\[
\vert x_{k}\vert^{p+d}\vert K(x)\vert\le\frac{1}{(2\pi)^{d}}\Big\Vert\frac{\partial^{p+d}\F K}{\partial u_{k}^{p+d}}\Big\Vert_{L^{1}(I_{1})}=\frac{1}{(2\pi)^{d}}\int_{I_{1}}\Big\vert\int\e^{\i\langle u,z\rangle}K(z)z_{k}^{p+d}\,\d z\Big\vert\,\d u\lesssim1,
\]
It follows from the equivalence of norms $\vert x\vert\lesssim\vert x\vert_{p+d}$
that
\begin{equation}
\vert x\vert^{p+d}\vert K(x)\vert\lesssim\vert x\vert_{p+d}^{p+d}\vert K(x)\vert=\sum_{k=1}^{d}\vert x_{k}\vert^{p+d}\vert K(x)\vert\lesssim1.\label{eq:kernelorder}
\end{equation}
Thus,
\[
\vert K_{h}(x^{\ast})\vert\le h^{-d}\sup_{\vert x\vert\ge\vert x^{\ast}\vert/h}\vert K(x)\vert\le h^{-d}\sup_{\vert x\vert\ge\vert x^{\ast}\vert/h}\frac{\vert x\vert^{p+d}}{\vert x^{\ast}/h\vert^{p+d}}\vert K(x)\vert\lesssim h^{p}\vert x^{\ast}\vert^{-p-d}
\]
and since $U$ is bounded away from $0$, this gives a uniform bound
in $x^{\ast}$ of the order $h^{s}$ as $p\ge s$.

\subsubsection{Controlling the bias\label{subsec:bias}}

For $x^{\ast}\in U$ and $h\vert x\vert<\zeta$, we use a multivariate
Taylor expansion  of $g\coloneqq\vert\cdot\vert^{2}\nu\in C^{s}(U_{\zeta})$
around $x^{\ast}$ to obtain
\[
g(x^{\ast}-hx)-g(x^{\ast})=\sum_{0<\vert\beta\vert_{1}<\lfloor s\rfloor}\frac{g^{(\beta)}(x^{\ast})}{\beta!}(-hx)^{\beta}+\sum_{\vert\beta\vert_{1}=\lfloor s\rfloor}\frac{g^{(\beta)}(x^{\ast}-\tau_{x^{\ast}-hx}hx)}{\beta!}(-hx)^{\beta},
\]
for some $\tau_{x^{\ast}-hx}\in[0,1]$. The order of the kernel and
the Hölder regularity of $g$ yield

\begin{align*}
\vert B^{\nu}(x^{\ast})\vert & =\big\vert\big(K_{h}\ast(\vert\cdot\vert^{2}\nu)-\vert\cdot\vert^{2}\nu\big)(x^{\ast})\big\vert\\
 & =\Big\vert\int\big(g(x^{\ast}-hx)-g(x^{\ast})\big)K(x)\,\d x\Big\vert\\
 & \le\Big\vert\int_{\vert x\vert\ge\zeta/h}\Big(g(x^{\ast}-hx)-\sum_{\vert\beta\vert_{1}\le\lfloor s\rfloor}\frac{g^{(\beta)}(x^{\ast})}{\beta!}(-hx)^{\beta}\Big)K(x)\,\d x\Big\vert\\
 & \qquad\qquad\!\!\!+\Big\vert\int_{\vert x\vert<\zeta/h}\sum_{\vert\beta\vert_{1}=\lfloor s\rfloor}\frac{(-hx)^{\beta}}{\beta!}\big(g^{(\beta)}(x^{\ast}-\tau_{x^{\ast}-hx}hx)-g^{(\beta)}(x^{\ast})\big)K(x)\,\d x\Big\vert\\
 & \lesssim\int_{\vert x\vert\ge\zeta/h}\vert g(x^{\ast}-hx)K(x)\vert\,\d x+\sum_{\vert\beta\vert_{1}\le\lfloor s\rfloor}\frac{\Vert g^{(\beta)}\Vert_{L^{\infty}(U)}}{\beta!}\frac{h^{s}}{\zeta^{s-\vert\beta\vert_{1}}}\int_{\vert x\vert\ge\zeta/h}\vert x\vert^{s}\vert K(x)\vert\,\d x\\
 & \qquad\qquad\!\!\!+\sum_{\vert\beta\vert_{1}=\lfloor s\rfloor}\frac{1}{\beta!}\int_{\vert x\vert<\zeta/h}\vert hx\vert^{\lfloor s\rfloor}\vert\tau_{x^{\ast}-hx}hx\vert^{s-\lfloor s\rfloor}\vert K(x)\vert\,\d x.
\end{align*}
The second and the third term are clearly of the order $h^{s}$. To
establish the same for the first term, we proceed slightly differently
for the two cases of ill-posedness.

In the severely ill-posed case, we separate the behavior of the small
and the large jumps. On the one hand, \ref{eq:kernelorder} yields
\begin{align*}
\int_{\substack{\vert x\vert\ge\zeta/h,\\
\vert x^{\ast}-hx\vert>1
}
}\vert g(x^{\ast}-hx)K(x)\vert\,\d x & \le\frac{h^{p+d}}{\zeta^{p+d}}\int_{\substack{\vert x\vert\ge\zeta/h,\\
\vert x^{\ast}-hx\vert>1
}
}\vert x^{\ast}-hx\vert^{2}\nu(x^{\ast}-hx)\vert x\vert^{p+d}\vert K(x)\vert\,\d x\\
 & \lesssim h^{s}\int_{\vert x\vert>1}\vert x\vert^{2}\nu(x)\,\d x.
\end{align*}

On the other hand, the assumption $\vert y\vert^{3-\eta}\nu(y)\le R$
for $\vert y\vert\le1$ and \ref{eq:kernelorder} gives 
\begin{align*}
\int_{\substack{\vert x\vert\ge\zeta/h,\\
\vert x^{\ast}-hx\vert\le1
}
}\vert g(x^{\ast}-hx)K(x)\vert\,\d x\lesssim\int_{\substack{\vert x\vert\ge\zeta/h,\\
\vert x^{\ast}-hx\vert\le1
}
}\frac{\vert x\vert^{-(p+d)}}{\vert x^{\ast}-hx\vert^{1-\eta}}\,\d x & \lesssim h^{p}\int_{\substack{\vert x^{\ast}-x\vert\ge\zeta,\\
\vert x\vert\le1
}
}\frac{\vert x^{\ast}-x\vert^{-(p+d)}}{\vert x\vert^{1-\eta}}\,\d x\\
 & \lesssim h^{s}\zeta^{-(p+d)}\int_{\vert x\vert\le1}\vert x\vert^{\eta-1}\,\d x.
\end{align*}

In the mildly ill-posed case, one uses $\Vert\vert x\vert\nu\Vert_{\infty}\le R$
and the triangle inequality to find 
\[
\int_{\vert x\vert\ge\zeta/h}\vert g(x^{\ast}-hx)K(x)\vert\,\d x\le R\int_{\vert x\vert\ge\zeta/h}\vert x^{\ast}-hx\vert\vert K(x)\vert\,\d x\lesssim h^{s}\frac{1+\zeta}{\zeta^{s}}\int\vert x\vert^{s}\vert K(x)\vert\,\d x.
\]

\subsubsection{Controlling the remainder term\label{subsec:remainder}}

To bound $R_{\tg,n}$ in \ref{eq:decomposition}, we first show that
with $a_{n}$ from \ref{lem:linearization} 
\[
g\coloneqq\F^{-1}\big[\F K_{h}\big(\hat{\lap\psi_{n}}-\lap\psi-\tg^{-1}\lap((\hat{\phi}_{\tg,n}-\phi_{\tg})/\phi_{\tg})\big)\big]=\mathcal{O}_{\P}\big(h^{-d}a_{n}\big).
\]
Let $\eps>0$. Owing to \ref{lem:linearization} we can choose $N,M>0$
(wlog. $M>\Vert K\Vert_{L^{1}}$) such that the probability of the
event $A_{n}\coloneqq\{\sup_{\vert u\vert_{\infty}\le h^{-1}}\vert\tilde g(u)\vert>a_{n}M\}$
with $\tilde g\coloneqq\hat{\lap\psi_{n}}-\lap\psi-\tg^{-1}\lap((\hat{\phi}_{\tg,n}-\phi_{\tg})/\phi_{\tg})$
is less than $\eps$ for $n>N$. Due to the support of $\F K$, we
have on $A_{n}^{\comp}$

\[
\vert g(x^{\ast})\vert=\frac{1}{(2\pi)^{d}}\Big\vert\int_{I_{h}}\e^{-\i\langle u,x^{\ast}\rangle}\F K(hu)\tilde g(u)\,\d u\Big\vert\le\frac{a_{n}M}{(2\pi)^{d}}\int_{I_{h}}\vert\F K(hu)\vert\,\d u\le h^{-d}a_{n}M\Vert K\Vert_{L^{1}}.
\]
For $M'\coloneqq M\Vert K\Vert_{L^{1}}$, we obtain
\[
\P\Big(\sup_{x^{\ast}\in U}\vert g(x^{\ast})\vert>h^{-d}a_{n}M'\Big)\le\eps,
\]
whereby the remainder term has the order proposed in \ref{eq:remainderorder}.

In the mildly ill-posed case, we have $\Vert\phi_{\tg}^{-1}\Vert_{L^{\infty}(I_{h})}\lesssim h^{-\tg\alpha}$
and \ref{eq:psigradwlog} implies $\big\Vert\vert\nabla\psi\vert\big\Vert_{L^{\infty}(I_{h})}\lesssim1$.
Thus, we have 
\[
\vert R_{\tg,n}\vert=\mathcal{O}_{\P}\big(n^{-1}\tg^{-1/2}(\log h^{-1})^{1+\chi}h^{-2\tg\alpha-d}\big).
\]

In the severely ill-posed case, $\Vert\phi_{\tg}^{-1}\Vert_{L^{\infty}(I_{h})}\lesssim\exp(r\tg h^{-\alpha})$
holds and \ref{eq:psigrad} implies $\big\Vert\vert\nabla\psi\vert\big\Vert_{L^{\infty}(I_{h})}\lesssim h^{-1}.$
Hence,
\[
\vert R_{\tg,n}\vert=\mathcal{O}_{\P}\big(n^{-1}\tg^{-1/2}(\log h^{-1})^{1+\chi}\big(h^{-d}+\tg h^{-d-1}+\tg^{3/2}h^{-d-2}\big)\exp(2r\tg h^{-\alpha})\big).
\]

In both cases, the remainder term is dominated by the linearized stochastic
error.

\medskip{}

This completes the proof of \ref{thm:main}.\hfill\qed

\subsection{Proof of \ref{prop:xs_nu}}

For the modified estimator, we have to replace $\tr(\Sigma)K_{h}$
with $\big(\tr(\Sigma)-\hat{\tr(\Sigma)}\big)K_{h}$ in the decomposition
\ref{eq:decomposition}. All other terms are treated as before. Since
we can bound $\vert K_{h}(x^{\ast})\vert$ by $h^{-d}$ uniformly
in $x^{\ast}$ and $\tg\le R$ is fixed, we only need to prove
\[
\big\vert\hat{\tr(\Sigma)}-\tr(\Sigma)\big\vert=\mathcal{O}_{\P}\big((\log n)^{-(s+d)/\alpha}\big).
\]
Similarly to \ref{eq:decomposition}, the error for estimating the
trace of $\Sigma$ can be decomposed into 
\begin{align*}
\tr(\Sigma)-\hat{\tr(\Sigma)} & =\int W_{h}(u)\big(\hat{\lap\psi_{n}}-\lap\psi\big)(u)\,\d u+\int W_{h}(u)\big(\lap\psi(u)+\tr(\Sigma)\big)\,\d u\\
 & =\underbrace{\int W_{h}(u)\tg^{-1}\lap\Big(\frac{\hat{\phi}_{\tg,n}(u)-\phi_{\tg}(u)}{\phi_{\tg}(u)}\Big)(u)\,\d u}_{\eqqcolon\tilde L_{\tg,n}^{\nu}}+\tilde R_{\tg,n}-\underbrace{\int W_{h}(u)\F[\vert x\vert^{2}\nu](u)\,\d u}_{\eqqcolon\tilde B_{h}^{\nu}}
\end{align*}
with the linearized stochastic error $\tilde L_{\tg,n}^{\nu}$, the
bias $\tilde B_{h}^{\nu}$ and a remainder term $\tilde R_{\tg,n}$.
Using the techniques from \ref{subsec:stochasticerror}, it is straightforward
to see

\begin{align*}
\E[\vert\tilde L_{\tg,n}^{\nu}\vert] & \lesssim\tg^{-1}n^{-1/2}\int\vert W_{h}(u)\vert\vert\phi_{\tg}^{-1}(u)\vert\big(\tg^{1/2}+\tg^{3/2}\vert\nabla\psi(u)\vert+\tg^{2}\vert\nabla\psi(u)\vert^{2}+\tg\big)\,\d u.\\
 & \lesssim\tg^{-1}n^{-1/2}\Vert\phi_{\tg}^{-1}\Vert_{L^{\infty}(I_{h})}\big(\tg^{1/2}+\tg^{3/2}h^{-1}+\tg^{2}h^{-2}\big)\int\vert W_{h}(u)\vert\,\d u\\
 & \lesssim n^{-1/2}\exp(r\tg h^{-\alpha})h^{-2}\int\vert W(u)\vert\,\d u,
\end{align*}
which is of the order $n^{-1/4}(\log n)^{2/\alpha}$ by our choice
of $h$ and will be dominated by the bias.

Using the Plancherel theorem in a similar fashion to \citet[Section 4.2.1]{belomestny2015},
we have for $g\coloneqq\vert\cdot\vert^{2}\nu$ and $\beta\in\N_{0}^{d}$
with $\beta_{l}=s\1_{\{l=k\}}$
\[
\vert\tilde B_{h}^{\nu}\vert\lesssim\Big\vert\int\F^{-1}[W_{h}(x)/x_{k}^{s}]\F^{-1}\big[x_{k}^{s}\F g(x)\big](u)\,\d u\Big\vert\lesssim\Vert g^{(\beta)}\Vert_{\infty}\big\Vert\F^{-1}[W_{h}(x)/x_{k}^{s}]\Vert_{L^{1}}.
\]
By substitution 
\[
\F^{-1}[W_{h}(x)/x_{k}^{s}](u)=\frac{h^{s}}{(2\pi)^{d}}\F^{-1}\big[W(x)/x_{k}^{s}\big](u/h)
\]
and therefore 
\[
\vert\tilde B_{h}^{\nu}\vert\lesssim h^{s}\Vert g^{(\beta)}\Vert_{\infty}\big\Vert\F^{-1}[W(x)/x_{k}^{s}](\cdot/h)\Vert_{L^{1}}\lesssim h^{(s+d)}\Vert g^{(\beta)}\Vert_{\infty}\big\Vert\F^{-1}[W(x)/x_{k}^{s}]\Vert_{L^{1}}\lesssim h^{(s+d)}.
\]

Together with \ref{lem:linearization}, we have 

\begin{align*}
\vert\tilde R_{\tg,n}\vert & \lesssim\sup_{\vert u\vert_{\infty}\le h^{-1}}\big\vert\hat{\lap\psi_{n}}(u)-\lap\psi(u)-\tg^{-1}\lap\big((\hat{\phi}_{\tg,n}-\phi_{\tg})/\phi_{\tg}\big)(u)\big\vert\int\vert W_{h}(u)\vert\,\d u\\
 & =\mathcal{O}_{\P}\big(n^{-1}(\log h^{-1})^{1+\chi}\Vert\phi_{\tg}^{-1}\Vert_{L^{\infty}(I_{h})}^{2}h^{-2}\big).
\end{align*}
which is dominated by the linearized stochastic error.\hfill\qed

\subsection{Proof of \ref{thm:misspec}}

The distributional analogon to \ref{eq:decomposition} is 
\begin{align*}
\int f(x)\vert x\vert^{2}\hat{\nu}_{h}(x)(\d x)-\int f(x)\vert x\vert^{2}\,\nu(\d x)\\
 & \negthickspace\negthickspace\negthickspace\negthickspace\negthickspace\negthickspace\negthickspace\negthickspace\negthickspace\negthickspace\negthickspace\negthickspace\negthickspace\negthickspace\negthickspace=\int f(x)\big(K_{h}\ast(\vert\cdot\vert^{2}\nu)\big)(x)\,\d x-\int f(y)\vert y\vert^{2}\nu(\d y)\\
 & \negthickspace\negthickspace\negthickspace\negthickspace\negthickspace\negthickspace\negthickspace\negthickspace\negthickspace\negthickspace\negthickspace\negthickspace\negthickspace\negthickspace\negthickspace\qquad+\int f(x)\F^{-1}\big[\F K_{h}(\lap\psi-\hat{\lap\psi}_{n})\big](x)\,\d x+\int f(x)\tr(\Sigma)K_{h}(x)\,\d x\\
 & \negthickspace\negthickspace\negthickspace\negthickspace\negthickspace\negthickspace\negthickspace\negthickspace\negthickspace\negthickspace\negthickspace\negthickspace\negthickspace\negthickspace\negthickspace=\int f(x)\big(K_{h}\ast(\vert\cdot\vert^{2}\nu)\big)(x)\,\d x-\int f(y)\vert y\vert^{2}\nu(\d y)\\
 & \negthickspace\negthickspace\negthickspace\negthickspace\negthickspace\negthickspace\negthickspace\negthickspace\negthickspace\negthickspace\negthickspace\negthickspace\negthickspace\negthickspace\negthickspace\qquad+\int f(x)L_{\tg,n}^{\nu}(x)\,\d x+\int f(x)R_{\tg,n}\,\d x+\int f(x)\tr(\Sigma)K_{h}(x)\,\d x
\end{align*}
with the same $L_{\tg,n}^{\nu}$ and $R_{\tg,n}$, for which we will
derive uniform upper bounds on $U$ which directly translate into
bounds when integrating against test functions due to their regularity.
For the integrated bias, we use Fubini's theorem to obtain
\begin{align*}
B_{I}^{\nu} & =\int f(x)\big(K_{h}\ast(\vert\cdot\vert^{2}\nu)\big)(x)\,\d x-\int f(y)\vert y\vert^{2}\nu(\d y)\\
 & =\int f(x)\Big(\int K_{h}(x-y)\vert y\vert^{2}\,\nu(\d y)\Big)\,\d x-\int f(y)\vert y\vert^{2}\,\nu(\d y)\\
 & =\int\big(\big(f\ast K_{h}(-\cdot)\big)(y)-f(y)\big)\vert y\vert^{2}\,\nu(\d y)\\
 & =\int\big(\big(K_{h}(-\cdot)\ast f\big)(y)-f(y)\big)\vert y\vert^{2}\,\nu(\d y).
\end{align*}
$\big(\big(K_{h}(-\cdot)\ast f\big)(y)-f(y)\big)$ is of the order
$(\vert y\vert\lor1)h^{\varrho}$ which follows from the arguments
in \ref{eq:decomposition} with $g=f$, $\varrho$ and $K_{h}(-\cdot)$
instead of $\vert\cdot\vert^{2}\nu$, $s$ and $K_{h}$, respectively.
Therefore,
\[
\vert B_{I}^{\nu}\vert\lesssim h^{\varrho}\int(\vert y\vert\lor1)\vert y\vert^{2}\,\nu(\d y)\lesssim h^{\varrho}\int\vert y\vert^{3}\,\nu(\d y)\lesssim h^{\varrho}.
\]

\medskip{}
A key tool to control linearized stochastic error $L_{\tg,n}^{\nu}$
in \ref{subsec:stochasticerror} was \ref{eq:moments}, which we can
still establish here by bounding the first four partial derivatives
of $\psi$ at the origin. Indeed, by \ref{eq:psilap} we have $\frac{\partial\psi}{\partial u_{k}}(0)=0$
and similarly

\begin{equation}
\Big\vert\frac{\partial\psi}{\partial u_{k}^{j}}(u)\Big\vert\le\tr(\Sigma)+\int\vert x\vert^{j}\,\nu(\d x)=\tr(\Sigma)+\sum_{l=1}^{2}\int\vert x_{l}\vert^{j}\,\nu_{l}(\d x_{l})\lesssim1,\quad j=2,3,4,\,k=1,\dots,d.\label{eq:psimoments}
\end{equation}
Hence, \ref{eq:moments} holds. Additionally, \ref{eq:psimoments}
implies $\vert\lap\psi(u)\vert\lesssim1$.

In the severely ill-posed case, we can still bound the gradient of
$\psi$ by 
\[
\vert\nabla\psi(u)\vert\lesssim\vert\Sigma u\vert+\int\vert x\vert\vert\e^{\i\langle u,x\rangle}-1\vert\,\nu(\d x)\le\vert u\vert+\int\vert\langle u,x\rangle\vert\vert x\vert\,\nu(\d x)\lesssim\vert u\vert,
\]
and then apply the arguments from \ref{eq:severely_linearized}. Hence,
the linearized stochastic error is of the same order as before.

In the severely ill-posed case, \ref{eq:psimoments} holds even for
$j=1$ and therefore $\vert\nabla\psi(u)\vert\lesssim1$. Continuing
from \ref{eq:moments} in the mildly ill-posed case requires the most
significant changes. \ref{eq:psigradwlog} now reads as
\[
\nabla\psi(u)=\big(\i\F[x^{(1)}\nu_{1}](u^{(1)}),\i\F[x^{(2)}\nu_{2}](u^{(2)})\big)^{\top}.
\]
and the main crux is that 
\begin{equation}
\int\e^{\i\langle u^{(k)},x^{(k)}\rangle}\vert x^{(k)}\vert^{j}\,\nu_{k}(\d x^{(k)}),\qquad j,k=1,2\label{eq:constant_arg}
\end{equation}
are constant in half of their arguments. Therefore, they cannot be
finitely integrable as functions on $\R^{d}$. In \ref{eq:diffstoch},
a way out is to consider
\begin{equation}
\big\Vert m_{\tg,h}\vert\nabla\psi\vert\big\Vert_{L^{1}}\le\sum_{k=1}^{2}\big\Vert m_{\tg,h}(u)\vert\F[x^{(k)}\nu_{k}](u^{(k)})\vert\big\Vert_{L^{1}}.\label{eq:m_aux}
\end{equation}
Then, we apply the Cauchy-Schwarz inequality and Plancherel's theorem
only on $L^{2}(\R^{d/2})$ to obtain
\begin{align}
\big\Vert m_{\tg,h}(u)\vert\F[x^{(1)}\nu_{1}](u^{(1)})\vert\big\Vert_{L^{1}}\nonumber \\
 & \negmedspace\negmedspace\negmedspace\negmedspace\negmedspace\negmedspace\negmedspace\negmedspace\negmedspace\negmedspace\negmedspace\negmedspace\negmedspace\negmedspace\negmedspace\negmedspace\negmedspace\negmedspace\negmedspace\negmedspace\negmedspace\negmedspace\negmedspace\negmedspace=\int\int\Big\vert\frac{\F K(hu)}{\phi_{\tg}(u)}\vert\F[x^{(1)}\nu_{1}](u^{(1)})\vert\Big\vert\,\d u^{(1)}\,\d u^{(2)}\\
 & \negmedspace\negmedspace\negmedspace\negmedspace\negmedspace\negmedspace\negmedspace\negmedspace\negmedspace\negmedspace\negmedspace\negmedspace\negmedspace\negmedspace\negmedspace\negmedspace\negmedspace\negmedspace\negmedspace\negmedspace\negmedspace\negmedspace\negmedspace\negmedspace\le\Vert\vert\F[x^{(1)}\nu_{1}]\vert\Vert_{L^{2}(\R^{d/2})}\int_{[-h^{-1},h^{-1}]^{d/2}}\Big(\int_{[-h^{-1},h^{-1}]^{d/2}}\big\vert\phi_{\tg}^{-1}(u)\big\vert^{2}\,\d u^{(1)}\Big)^{1/2}\,\d u^{(2)}\nonumber \\
 & \negmedspace\negmedspace\negmedspace\negmedspace\negmedspace\negmedspace\negmedspace\negmedspace\negmedspace\negmedspace\negmedspace\negmedspace\negmedspace\negmedspace\negmedspace\negmedspace\negmedspace\negmedspace\negmedspace\negmedspace\negmedspace\negmedspace\negmedspace\negmedspace\lesssim h^{-\tg\alpha-3d/4}.\label{eq:m_aux_grad}
\end{align}
Analogously, the second summand in \ref{eq:m_aux} has the same order.
As a direct consequence,
\[
\big\Vert m_{\tg,h}\vert\nabla\psi\vert^{2}\big\Vert_{L^{1}}\le\sum_{k=1}^{2}\Vert\vert x^{(k)}\vert\nu_{k}\Vert_{L^{1}}\big\Vert m_{\tg,h}\vert\nabla\psi\vert\big\Vert_{L^{1}}\lesssim h^{-\tg\alpha-3d/4}
\]
and similarly to \ref{eq:m_aux_grad,eq:m_aux} the same holds for
$\Vert m_{\tg,h}\lap\psi\Vert_{L^{1}}$. Recalling \ref{eq:diffstoch},
we have 
\begin{align*}
\E\Big[\sup_{x^{\ast}\in U}\big\vert M_{\tg,n}^{\nu}(x^{\ast})-L_{\tg,n}^{\nu}(x^{\ast})\big\vert\Big] & \le(2\pi)^{-d}\big(2\E\big[\big\Vert m_{\tg,h}\nabla(\hat{\phi}_{\tg,n}-\phi_{\tg})\cdot\nabla\psi\big\Vert_{L^{1}}\big]\\
 & \qquad\qquad\qquad+\E\big[\big\Vert m_{\tg,h}(\hat{\phi}_{\tg,n}-\phi_{\tg})\big(\lap\psi-\tg(\nabla\psi)^{2}\big)\big\Vert_{L^{1}}\big]\big)\\
 & \lesssim n^{-1/2}h^{-\tg\alpha-3d/4}.
\end{align*}
Note that we pay for the dependence struture with an additional $h^{-d/4}$
compared to \ref{eq:plancherel}. The same happens when applying Bernstein's
inequality to obtain the following adaptation of \ref{lem:bernsteinapplication}.
\begin{lem}
\label{lem:bernstein_misspec}Let $\alpha,R,\zeta>0,m>4,U\subseteq\R^{d}$
and $x\in\R^{d}$, let the kernel have product structure and satisfy
\ref{eq:kernel} for $p\ge1$. If $(\Sigma,\gamma,\nu)\in\tilde{\mathcal{D}}(\alpha,m,R)$,
then there exists some constant $c>0$ depending only on $R,\alpha$
and $d$ such that for any $\kappa_{0}>0$ and any $n\in\N,0<\tg\le R,h\in(0,1)$
\[
\P\big(\vert M_{\tg,n}^{\nu}(x)\vert\ge\kappa_{0}T^{-1/2}h^{-\tg\alpha-3d/4}\big)\le2\exp\Big(-\frac{c\kappa_{0}^{2}}{(1+\vert x\vert^{3})(1+\kappa_{0}(Th^{d/2})^{-1/2})}\Big).
\]
\end{lem}

Carrying out the discretization argument from before, the linearized
stochastic error in the mildly ill-posed case is of the order 
\[
\vert L_{\tg,n}^{\nu}\vert=\mathcal{O}_{\P}\Big(\Big(\frac{\log T}{T}\Big)^{1/2}h^{-\tg\alpha-3d/4}\Big).
\]
\medskip{}
The term $\tr(\Sigma)K_{h}$ is treated as in \ref{subsec:volatility}
just with $\varrho$ instead of $s$. No changes are necessary to
treat the remainder term compared to \ref{subsec:remainder}. This
is because when treating the linearized stochstic error, we already
showed that still $\vert\nabla\psi(u)\vert,\vert\lap\psi(u)\vert\lesssim1$
in the mildly ill-posed case and $\vert\nabla\psi(u)\vert\lesssim\vert u\vert,\vert\lap\psi(u)\vert\lesssim1$
in the severely ill-posed case. This concludes the proof of \ref{thm:misspec}.\hfill\qed

\subsection{Remaining proofs}

\subsubsection{Proof of \ref{lem:characteristicinequality}}

The proof uses empirical process theory and is a combination of \citet{Kappus2010}
and \citet{belomestny2018}.

To simplify the notation, write
\[
C_{\rho,n}^{\beta}(u)\coloneqq n^{-1/2}\rho^{-(\vert\beta\vert_{1}\land1)/2}\sum_{k=1}^{n}\frac{\partial^{\beta}}{\partial u^{\beta}}\big(\e^{\i\langle u,X_{k}\rangle}-\E[\e^{\i\langle u,X_{k}\rangle}]\big)
\]
so that the assertion reads
\[
\sup_{\substack{n\ge1,\\
0<\rho\le R
}
}\E\big[\big\Vert w(u)C_{\rho,n}^{\beta}(u)\big\Vert_{\infty}\big]<\infty.
\]
We decompose $C_{\rho,n}^{\beta}$ into its real and its imaginary
part to obtain
\[
\E\big[\big\Vert w(u)C_{\rho,n}^{\beta}(u)\big\Vert_{\infty}\big]\le\E\big[\big\Vert w(u)\Re\big(C_{\rho,n}^{\beta}(u)\big)\big\Vert_{\infty}\big]+\E\big[\big\Vert w(u)\Im\big(C_{\rho,n}^{\beta}(u)\big)\big\Vert_{\infty}\big].
\]
As both parts can be treated analogously, we focus on the real part.
To this end, introduce the class of 
\[
\mathcal{G}_{\rho,\beta}\coloneqq\{g_{u}:u\in\R^{d}\}\qquad\text{where}\qquad g_{u}\colon\R^{d}\to\R,\qquad x\mapsto w(u)\rho^{-(\vert\beta\vert_{1}\land1)}\frac{\partial^{\beta}}{\partial u^{\beta}}\cos(\langle u,x\rangle).
\]
Since $G=\rho^{-(\vert\beta\vert_{1}\land1)/2}\vert\cdot\vert^{\beta}$
is an envelope function for $\mathcal{G}_{\rho,\beta}$, \citet[Corollary 19.35]{Vaart1998}
yields
\[
\E\big[\big\Vert w(u)\Re\big(C_{\rho,n}^{\beta}(u)\big)\big\Vert_{\infty}\big]\lesssim J_{[]}\big(\E[G(X_{1})^{2}]^{1/2},\mathcal{G}_{\rho,\beta}\big)\coloneqq\int_{0}^{\E[G(X_{1})^{2}]^{1/2}}\sqrt{\log N_{[]}(\eps,\mathcal{G}_{\rho,\beta})}\,\d\eps,
\]
where $N_{[]}(\eps,\mathcal{G}_{\rho,\beta})$ is the minimal number
of $\eps$-brackets (with respect to the distribution of $X_{1}$)
needed to cover $\mathcal{G}_{\rho,\beta}$. 

Since $\vert g_{u}(x)\vert\le w(u)\rho^{-(\vert\beta\vert_{1}\land1)/2}\vert x\vert^{\beta}$,
the set $\{g_{u}:\vert u\vert>B\}$ is covered by the bracket
\begin{gather*}
[g_{0}^{-},g_{0}^{+}]\coloneqq\{g\colon\R^{d}\to\R\mid g_{0}^{-}(x)\le g(x)\le g_{0}^{+}(x)\,\forall x\in\R^{d}\}\qquad\text{for}\\
g_{0}^{\pm}\coloneqq\pm\eps\rho^{-(\vert\beta\vert_{1}\land1)/2}\vert\cdot\vert^{\beta}\qquad\text{and}\qquad B\coloneqq B(\eps)\coloneqq\inf\big\{ b>0:\sup_{\vert u\vert\ge b}w(u)\le\eps\big\}.
\end{gather*}
To cover $\{g_{u}:\vert u\vert\le B\}$, we use for some grid $(u_{\rho,j})_{j\ge1}\subseteq\R^{d}$
the functions

\[
g_{\rho,j}^{\pm}\coloneqq\rho^{-(\vert\beta\vert_{1}\land1)/2}\big(w(u_{\rho,j})\frac{\partial^{\beta}}{\partial u_{\rho,j}^{\beta}}\cos(\langle u_{\rho,j},\cdot\rangle)\pm\eps\vert\cdot\vert^{\beta}\big)\1_{\{\vert\cdot\vert\le M\}}\pm\rho^{-(\vert\beta\vert_{1}\land1)/2}\vert\cdot\vert^{\beta}\1_{\{\vert\cdot\vert>M\}}.
\]
where $M\coloneqq\inf\{m:\rho^{-(\vert\beta\vert_{1}\land1)}\E[\vert X_{1}\vert^{2\beta}\1_{\{\vert X_{1}\vert>m\}}]\le\eps^{2}\}$.
Owing to $\E[\vert X_{1}\vert^{2\beta}]\lesssim\rho^{\vert\beta\vert_{1}\land1}$,
we have 
\[
\E\big[\vert g_{j}^{+}(X_{1})-g_{j}^{-}(X_{1})\vert^{2}\big]\le4\eps^{2}\big(\rho^{-(\vert\beta\vert_{1}\land1)}\E[\vert X_{1}\vert^{2\beta}]+1\big)\le c\eps^{2}
\]
for some $c>0$. Denote by $C$ the Lipschitz constant of $w$ and
use the triangle inequality to see
\[
\Big\vert w(u)\frac{\partial^{\beta}}{\partial u^{\beta}}\cos(\langle u,x\rangle)-w(u_{j})\frac{\partial^{\beta}}{\partial u_{\rho,j}^{\beta}}\cos(\langle u_{\rho,j},x\rangle)\Big\vert\le\vert x\vert^{\beta}(C+\vert x\vert)\vert u-u_{\rho,j}\vert.
\]
Thus, $g_{u}\in[g_{j}^{-},g_{j}^{+}]$ as soon as $(C+M)\vert u-u_{\rho,j}\vert\le\eps$.
It takes at most $(\lceil B/\eps_{0}\rceil)^{d}$ $\ell^{2}$-balls
of radius $d^{1/2}\eps_{0}$ to cover the $\ell^{2}$-ball of radius
$B$ around $0$. For $\eps_{0}=\eps d^{-1/2}/(C+M)$, denote their
centers by $(u_{\rho,j})_{j}$. To translate this into a cover of
$\{g_{u}:\vert u\vert\le B\}$, we fix some $g_{u}$ with $\vert u\vert\le B$.
By construction, we can pick $j$ such that $\vert u-u_{\rho,j}\vert\le d^{1/2}\eps_{0}=\eps/(C+M)$.
The previous calculations show that $[g_{j}^{-},g_{j}^{+}]$ is a
$c^{1/2}\eps$-bracket containing $g_{u}$ and therefore 
\[
N_{[]}(\eps,\mathcal{G}_{\rho,\beta})\le\big(\big\lceil\eps^{-1}(cd)^{1/2}B(C+M)\big\rceil\big)^{d}+1.
\]

It is straightforward to see $B\le\exp(\eps^{-2/(1+\chi)}$). Further,
$m=(\eps^{-2}\rho^{-(\vert\beta\vert_{1}\land1)}\E[\vert X_{1}\vert^{2\beta}\vert X_{1}\vert^{\tau}])^{1/\tau}$
is sufficient for
\[
\rho^{-(\vert\beta\vert_{1}\land1)}\E[\vert X_{1}\vert^{2\beta}\1_{\{\vert X_{1}\vert>m\}}]\le m^{-\tau}\E[\vert X_{1}\vert^{2\beta}\vert X_{1}\vert^{\tau}]\le\eps^{2}
\]
and thus $M\le(\eps^{-2}\rho^{-(\vert\beta\vert_{1}\land1)}\E[\vert X_{1}\vert^{2\beta}\vert X_{1}\vert^{\tau}])^{1/\tau}\le(\eps^{-2}dc')^{1/\tau}$
for some $c'>0$. Hence,
\[
\log N_{[]}(\eps,\mathcal{G}_{\rho,\beta})\lesssim1+\log(\eps^{-2/\tau-1})+\eps^{-2/(1+\tau)}\lesssim1+\eps^{-2/(1+\tau)}
\]
implying
\[
J_{[]}\big(\E[G(X_{1})^{2}]^{1/2},\mathcal{G}_{\rho,\beta}\big)=\int_{0}^{(\rho^{-(\vert\beta\vert_{1}\land1)}\E[\vert X_{1}\vert^{2\beta}])^{1/2}}\sqrt{\log N_{[]}(\eps,\mathcal{G}_{\rho,\beta})}\,\d\eps<\infty.\tag*{{\qed}}
\]

\subsubsection{Proof of \ref{lem:linearization}}

Setting $g(y)=\log(1+y)$ (ie. $g'(y)=(1+y)^{-1},g''(y)=-(1+y)^{-2})$
and $\xi=(\hat{\phi}_{\tg,n}-\phi_{\tg})/\phi_{\tg}$, we use
\[
\nabla(g\circ\xi)(u)=g'(\xi(u))\nabla\xi(u),\qquad\lap(g\circ\xi)(u)=g''(\xi(u))(\nabla\xi)^{2}(u)+g'(\xi(u))\lap\xi(u)
\]
and $\vert(\nabla\xi)^{2}(u)\vert\le\vert\nabla\xi(u)\vert^{2}$ to
obtain for $\vert\xi(u)\vert\le1/2$ that
\begin{equation}
\big\vert\lap(g\circ\xi)(u)-\lap\xi(u)\big\vert\le\vert g''(\xi(u))\vert\vert\nabla\xi(u)\vert^{2}+\vert g'(\xi(u))-1\vert\vert\lap\xi(u)\vert\lesssim\vert\nabla\xi(u)\vert^{2}+\vert\xi(u)\vert\vert\lap\xi(u)\vert,\label{eq:prelinear}
\end{equation}
because
\[
\vert g'(y)-1\vert\le2\vert y\vert\qquad\text{and}\qquad\vert g''(y)\vert\le4\qquad\forall y\in\C:\vert y\vert\le1/2.
\]
The latter statement holds, since $1/2\le\vert1+y\vert$. For the
former statement, consider the expansion 
\[
g'(y)=\frac{1}{1+y}=\sum_{k=0}^{\infty}(-y)^{k}\qquad\forall y\in\C:\vert y\vert\le1/2
\]
to see 
\[
\vert g'(y)-1\vert=\Big\vert\sum_{k=1}^{\infty}(-y)^{k}\Big\vert=\Big\vert-y\sum_{k=0}^{\infty}(-y)^{k}\Big\vert=\Big\vert\frac{y}{1+y}\Big\vert\le2\vert y\vert.
\]

Noting $\hat{\lap\psi_{n}}-\lap\psi=\tg^{-1}\lap\log(\hat{\phi}_{\tg,n}/\phi_{\tg})$,
\ref{eq:prelinear} implies on the event $\Omega_{n}\coloneqq\big\{\sup_{\vert u\vert_{\infty}\le h^{-1}}\vert\xi(u)\vert\le1/2\big\}$
\[
\sup_{\vert u\vert_{\infty}\le h^{-1}}\big\vert\tg(\hat{\lap\psi_{n}}-\lap\psi)(u)-\lap((\hat{\phi}_{\tg,n}-\phi_{\tg})/\phi_{\tg})(u)\big\vert\lesssim\Vert\vert\nabla\xi\vert\Vert_{L^{\infty}(I_{h})}^{2}+\Vert\xi\Vert_{L^{\infty}(I_{h})}\Vert\lap\xi\Vert_{L^{\infty}(I_{h})}.
\]

To control the $\xi$-terms, we invoke \ref{lem:characteristicinequality}
applied to the increments of the Lévy process with $\rho=\tg$ after
verifying that the moments are of the appropriate order. Owing to
the equivalence of norms, it is sufficient to show that with $\tau=m-4>0$
\begin{equation}
\E[\vert Y_{1,k}\vert^{2l+\tau}]\lesssim\tg^{l\land1}\qquad\text{and}\qquad\E[\vert Y_{1,k}\vert^{2l}]\lesssim\tg^{l\wedge1},\qquad k=1,\dots,d,\,l=0,1,2,\label{eq:momentorder}
\end{equation}
where $Y_{1,k}$ is the $k$-th entry of $Y_{1}$ and thus an increment
with time difference $\tg$ based on the Lévy process $(L_{t,k})_{t\ge0}$
with Lévy measure $\nu_{k}$. For $l=1,2$, it follows from \citet[Theorem 1.1]{FigueroaLopez2008}
that
\[
\lim_{\tg\searrow0}\tg^{-1}\E[\vert Y_{1,k}\vert^{2l+\tau}]=\lim_{\tg\searrow0}\tg^{-1}\E[\vert L_{\tg,k}\vert^{2l+\tau}]=\int\vert x_{k}\vert^{2l+\tau}\,\nu_{k}(\d x_{k})\le\int\vert x\vert^{2l+\tau}\,\nu(\d x)\lesssim R.
\]
For $l=0$, $\E[\vert Y_{1,k}\vert^{\tau}]\lesssim\E[\vert Y_{1,k}\vert^{m}]\lesssim1$
holds by our moment assumptions. The second condition in \ref{eq:momentorder}
was already checked at the beginning of \ref{subsec:stochasticerror}. 

Therefore, $\vert\lap\psi(u)\vert\lesssim1$ yields
\begin{align}
\Vert\xi\Vert_{L^{\infty}(I_{h})} & =\mathcal{O}_{\P}\big(n^{-1/2}(\log h^{-1})^{(1+\chi)/2}\Vert\phi_{\tg}^{-1}\Vert_{L^{\infty}(I_{h})}\big),\label{eq:etabound}\\
\Vert\vert\nabla\xi\vert\Vert_{L^{\infty}(I_{h})}^{2} & =\mathcal{O}_{\P}\big(n^{-1}(\log h^{-1})^{1+\chi}\Vert\phi_{\tg}^{-1}\Vert_{L^{\infty}(I_{h})}^{2}\big(\tg+\tg^{2}\Vert\vert\nabla\psi\vert\Vert_{L^{\infty}(I_{h})}^{2}\big)\big),\nonumber \\
\Vert\lap\xi\Vert_{L^{\infty}(I_{h})} & =\mathcal{O}_{\P}\big(n^{-1/2}(\log h^{-1})^{(1+\chi)/2}\Vert\phi_{\tg}^{-1}\Vert_{L^{\infty}(I_{h})}\big(\tg^{1/2}+\tg^{3/2}\big\Vert\vert\nabla\psi\vert\big\Vert_{L^{\infty}(I_{h})}+\tg^{2}\big\Vert\vert\nabla\psi\vert\big\Vert_{L^{\infty}(I_{h})}^{2}\big)\big).\nonumber 
\end{align}
Combining \ref{eq:etabound} with $n^{-1/2}(\log h^{-1})^{(1+\chi)/2}\Vert\phi_{\tg}^{-1}\Vert_{L^{\infty}(I_{h})}\to0$
gives $\P(\Omega_{n})\to1$, completing the proof.\hfill\qed

\subsubsection{Proof of \ref{lem:bernsteinapplication}}

For fixed $x\in\R^{d}$, we want to apply Bernstein's inequality to
\[
M_{\tg,n}^{\nu}(x)=-\sum_{l=1}^{n}\big(\xi_{l}-\E[\xi_{l}]\big)\qquad\text{with}\qquad\xi_{l}\coloneqq T^{-1}\F^{-1}\big[m_{\tg,h}(u)\vert Y_{l}\vert^{2}\e^{\i\langle u,Y_{l}\rangle}\big](x).
\]

Similar arguments to \ref{eq:m_order} reveal $\big\Vert m_{\tg,h}(u)\big\Vert_{L^{1}}\lesssim h^{-\tg\alpha-d}$,
and with the quotient rule one finds the same order for $\big\Vert\lap m_{\tg,h}(u)\big\Vert_{L^{1}}$
and $\big\Vert\vert\nabla m_{\tg,h}\vert\big\Vert_{L^{1}}$ paving
a deterministic bound of $\xi_{l}$ via
\begin{align*}
\vert\xi_{l}\vert & =T^{-1}\vert Y_{l}\vert^{2}\big\vert\F^{-1}\big[m_{\tg,h}(u)\e^{-\i\langle u,x\rangle}\big](-Y_{l})\big\vert\\
 & =T^{-1}\big\vert\F^{-1}\big[\lap\big(m_{\tg,h}(u)\e^{-\i\langle u,x\rangle}\big)\big](-Y_{l})\big\vert\\
 & \le T^{-1}\big\Vert\lap\big(m_{\tg,h}(u)\e^{-\i\langle u,x\rangle}\big)\big\Vert_{L^{1}}\\
 & \le T^{-1}\big(\big\Vert\lap m_{\tg,h}(u)\big\Vert_{L^{1}}+2\vert x\vert_{1}\big\Vert\vert\nabla m_{\tg,h}\vert\big\Vert_{L^{1}}+\vert x\vert^{2}\Vert m_{\tg,h}\Vert_{L^{1}}\big)\\
 & \lesssim T^{-1}(1+\vert x\vert^{2})h^{-\tg\alpha-d}.
\end{align*}

To bound the variance of $\xi_{l}$, note that for the distribution
$\P_{\tg}$ of $Y_{1}$, we have 
\[
\F[\i z_{k}\P_{\tg}]=\frac{\partial\phi_{\tg}}{\partial u_{k}}=\tg\phi_{\tg}\frac{\partial\psi}{\partial u_{k}}=\tg\F[\i z_{k}\nu]\phi_{\tg}=\F[\F^{-1}[\tg\F[\i z_{k}\nu]\phi_{\tg}]]=\F[(\tg\i z_{k}\nu)\ast\P_{\tg}]
\]
and therefore $z_{k}\P_{\tg}=\tg\mu\ast\P_{\tg}$ with $\mu(\d z)=z_{k}\nu(z)\d z$.
It follows that
\[
\int g(z)\vert z_{k}\vert\,\P_{\tg}(\d z)\le\tg\Vert z_{k}\nu\Vert_{\infty}\Vert g\Vert_{L^{1}},\qquad\forall g\in L^{1}(\R^{d}).
\]
Again, using similar arguments to \ref{eq:m_order} and the quotient
rule, we also have $\big\Vert\lap m_{\tg,h}(u)\big\Vert_{L^{2}},\big\Vert\vert\nabla m_{\tg,h}\vert\big\Vert_{L^{2}}\lesssim h^{-\tg\alpha-d/2}$.
Thus, the Cauchy-Schwarz inequality and the Plancherel theorem imply
\begin{align*}
\Var(\xi_{l})\le\E[\vert\xi_{l}\vert^{2}] & =T^{-2}\E\Big[\vert Y_{l}\vert^{4}\big\vert\F^{-1}\big[m_{\tg,h}(u)\e^{-\i\langle u,x\rangle}\big](-Y_{l})\big\vert^{2}\Big]\\
 & \lesssim T^{-2}\sum_{k=1}^{d}\int\vert y\vert^{3}\big\vert\F^{-1}\big[m_{\tg,h}(u)\e^{-\i\langle u,x\rangle}\big](-y)\big\vert^{2}\vert y_{k}\vert\,\P_{\tg}(\d y)\\
 & \le n^{-2}\tg^{-1}\sum_{k=1}^{d}\Vert z_{k}\nu\Vert_{\infty}\int\vert y\vert^{3}\big\vert\F^{-1}\big[m_{\tg,h}(u)\e^{-\i\langle u,x\rangle}\big](y)\big\vert^{2}\,\d y\\
 & \lesssim n^{-2}\tg^{-1}\big\Vert\vert y\vert^{2}\F^{-1}\big[m_{\tg,h}(u)\e^{-\i\langle u,x\rangle}\big](y)\big\Vert_{L^{2}}\big\Vert\vert y\vert_{1}\F^{-1}\big[m_{\tg,h}(u)\e^{-\i\langle u,x\rangle}\big](y)\big\Vert_{L^{2}}\\
 & \lesssim n^{-2}\tg^{-1}\Big(\sum_{k=1}^{d}\Big\Vert\frac{\partial^{2}}{\partial u_{k}^{2}}\big(m_{\tg,h}(u)\e^{-\i\langle u,x\rangle}\big)\Big\Vert_{L^{2}}\Big)\Big(\sum_{k=1}^{d}\Big\Vert\frac{\partial}{\partial u_{k}}\big(m_{\tg,h}(u)\e^{-\i\langle u,x\rangle}\big)\Big\Vert_{L^{2}}\Big)\\
 & \lesssim n^{-2}\tg^{-1}h^{-2\tg\alpha-d}(1+\vert x\vert^{3}).
\end{align*}

Now, Bernstein's inequality, e.g. \citet[Lemma 19.32]{Vaart1998}
yields for a constant $c'>0$ and any $\kappa>0$ that
\[
\P\big(\vert M_{\tg,n}^{\nu}(x)\vert\ge\kappa\big)\le2\exp\Big(-\frac{Tc'\kappa^{2}}{h^{-2\tg\alpha-d}(1+\vert x\vert^{3})+\kappa(1+\vert x\vert^{2})h^{-\tg\alpha-d}}\Big),
\]
which reads as the assertion if we choose $\kappa=\kappa_{0}T^{-1/2}h^{-\tg\alpha-d/2}$
for any $\kappa_{0}>0$ and set $c=c'/2$.\hfill\qed

\subsubsection{Proof of \ref{lem:bernstein_misspec}}

Fix $x=(x^{(1)},x^{(2)})$ for $x^{(1)},x^{(2)}\in\R^{d/2}$ and analogously
split $Y_{l}$ into its first and last $d/2$ entries $Y_{l}^{(1)}$
and $Y_{l}^{(2)}$ with characteristic functions $\phi_{\tg,1}$ and
$\phi_{\tg,2}$, respectively. Due to the product kernel, we obtain
\[
\xi_{l}=T^{-1}\vert Y_{l}\vert^{2}\big\vert\F^{-1}\big[m_{\tg,h}(u)\e^{-\i\langle u,x\rangle}\big](-Y_{l})\big\vert=T^{-1}(A_{1}B_{2}+A_{2}B_{1})
\]
with 
\begin{gather*}
A_{k}\coloneqq\vert Y_{l}^{(k)}\vert^{2}\F^{-1}[m_{\tg,h}^{(k)}(u^{(k)})\e^{\i\langle u^{(k)},Y_{l}^{(k)}\rangle}](x^{(k)}),\qquad B_{k}\coloneqq\F^{-1}[m_{\tg,h}^{(k)}(u^{(k)})\e^{\i\langle u^{(k)},Y_{l}^{(k)}\rangle}](x^{(k)}),\qquad\text{and}\\
m_{\tg,h}^{(k)}(u^{(k)})\coloneqq\phi_{\tg,k}^{-1}(u^{(k)})\prod_{j=1+(k-1)d/2}^{(k+1)d/2}\F K^{j}(hu_{j}^{(k)}),\qquad k=1,2.
\end{gather*}
$A_{1}$ and $A_{2}$ are the same terms that appeared in the proof
of \ref{lem:bernsteinapplication} just with half the dimension and
therefore
\begin{align*}
\vert A_{k}\vert & \lesssim\Vert\phi_{\tg,k}^{-1}\Vert_{L^{\infty}([-h^{-1},h^{-1}]^{d/2}])}h^{-d/2}(1+\vert x^{(k)}\vert^{2}),\\
\E[\vert A_{k}\vert^{2}] & \lesssim\tg\Vert\phi_{\tg,k}^{-1}\Vert_{L^{\infty}([-h^{-1},h^{-1}]^{d/2}])}^{2}h^{-d/2}(1+\vert x^{(k)}\vert^{3}).
\end{align*}
In a similar vain, $m_{\tg,h}^{(1)}$ and $m_{\tg,h}^{(2)}$ can be
treated be treated like $m_{\tg,h}$ with half the dimension leading
to
\[
\vert B_{k}\vert\lesssim\Vert m_{\tg,h}^{(k)}\Vert_{L^{1}(\R^{d/2})}\lesssim\Vert\phi_{\tg,k}^{-1}\Vert_{L^{\infty}([-h^{-1},h^{-1}]^{d/2}])}h^{-d/2}.
\]

Note that
\[
\prod_{k=1}^{2}\Vert\phi_{\tg,k}^{-1}\Vert_{L^{\infty}([-h^{-1},h^{-1}]^{d/2}])}=\Vert\phi_{\tg}^{-1}\Vert_{L^{\infty}(I_{h})}.
\]
Together, we have the deterministic bound $\vert\xi_{l}\vert\lesssim T^{-1}h^{-\tg\alpha-d}(1+\vert x\vert^{2})$.
Further, since $A_{1}$ and $B_{2}$ as well as $A_{2}$ and $B_{1}$
are independent, we obtain
\begin{align*}
\Var(\xi_{l})\lesssim T^{-2}(\Var(A_{1}B_{2})+\Var(A_{2}B_{1})) & \le T^{-2}\big(\E[\vert A_{1}\vert^{2}]\E[\vert B_{2}\vert^{2}]+\E[\vert A_{2}\vert^{2}]\E[\vert B_{1}\vert^{2}]\big)\\
 & \lesssim n^{-2}\tg^{-1}h^{-2\tg\alpha-3d/2}(1+\vert x\vert^{3}).
\end{align*}

Overall, Bernstein's inequality gives for a constant $c'>0$ and any
$\kappa>0$ that
\[
\P\big(\vert M_{\tg,n}^{\nu}(x)\vert\ge\kappa\big)\le2\exp\Big(-\frac{Tc'\kappa^{2}}{h^{-2\tg\alpha-d}(1+\vert x\vert^{3})+\kappa(1+\vert x\vert^{2})h^{-\tg\alpha-d}}\Big).
\]
The assertion follows by choosing $\kappa=\kappa_{0}T^{-1/2}h^{-\tg\alpha-3d/4}$
for any $\kappa_{0}>0$ and setting $c=c'/2$.\hfill\qed

\bibliographystyle{apalike2}
\bibliography{sources}

\end{document}